\documentclass[a4paper, 11pt]{article}
\usepackage{amsmath, amssymb, amsthm}
\usepackage{graphicx}
\usepackage{tikz}
\usetikzlibrary{calc, positioning, arrows.meta, decorations.markings}
\usepackage{float}
\usepackage{placeins}
\usepackage[skip=6pt,font=small]{caption}
\usepackage[colorlinks=true, linkcolor=blue, citecolor=blue, urlcolor=blue]{hyperref}

\newtheorem{theorem}{Theorem}[section]
\newtheorem{lemma}[theorem]{Lemma}
\newtheorem{corollary}[theorem]{Corollary}
\newtheorem{proposition}[theorem]{Proposition}

\theoremstyle{remark}
\newtheorem{remark}[theorem]{Remark}

\newcommand{\E}{\mathbb{E}}
\newcommand{\Prob}{\mathbb{P}}

\title{The number of cycles in a random 2-factor\\ of a random regular graph}
\author{%
Katerina Adler\\
{\normalsize Department of Mathematics, Bar-Ilan University, Ramat Gan, Israel}
\and
Simi Haber\\
{\normalsize Department of Mathematics, Bar-Ilan University, Ramat Gan, Israel}%
}
\date{August 2026}

\begin{document}
\maketitle

\begin{abstract}
We study the number of cycles $X$ in a uniformly random 2-factor $F$ of a random $d$-regular graph on $n$ vertices (in the pairing model) for any fixed $d\ge 3$. When the pairing $P$ and the 2-factor $F\subseteq P$ are chosen uniformly at random together, the probability generating function of $X$ is $\Gamma(n+\tfrac u2)\Gamma(\tfrac12)/\bigl(\Gamma(n+\tfrac12)\Gamma(\tfrac u2)\bigr)$, which does not depend on $d$. Consequently, $X$ has mean $\tfrac12\log n+O(1)$ and variance $\tfrac12\log n+O(1)$, satisfies a central limit theorem, and equals $1$ (so that $F$ is a Hamilton cycle) with probability $\Theta(n^{-1/2})$. Using an exact formula for the second moment of the cycle generating polynomial together with small subgraph conditioning, we transfer these results to a uniformly random 2-factor of a fixed random regular graph.
\end{abstract}


\section{Introduction}
A \emph{2-factor} of a graph is a spanning 2-regular subgraph: a disjoint union of cycles covering every vertex. In random regular graphs, 2-factors play a special algorithmic role as ``almost Hamilton cycles.'' The classical route of Frieze, Jerrum, Molloy, Robinson, and Wormald \cite{frieze1996} to generating Hamilton cycles in a random $d$-regular graph $G(n,d)$ is to sample near-uniform 2-factors repeatedly (via Tutte's reduction of 2-factors to perfect matchings \cite{tutte1954} and a matching algorithm) until a Hamiltonian one appears; their counting argument shows that the number of 2-factors exceeds the number of Hamilton cycles by at most a polynomial factor with high probability, so polynomially many resamplings suffice. Structurally, random regular graphs contain 2-factors with high probability \cite{robalewska1996}, admit 1-factorizations with high probability \cite{molloy1997}, and in fact decompose with high probability into Hamilton cycles and at most one perfect matching \cite{kim_wormald2001}; the contiguity toolbox behind such results is surveyed in \cite{greenhill2002,janson1995,wormald1999}. Recently, the cycle structure of random cycle-factors has attracted renewed attention: Christoph, Dragani\'c, Gir\~ao, Hurley, Michel, and M\"uyesser \cite{christoph2025} proved that a uniformly random cycle-factor of an \emph{arbitrary} $d$-regular digraph has $\Theta((n\log d)/d)$ cycles in expectation, and Kahn, Lubetzky, and Wormald \cite{kahn2017} studied $k$-cycle factors of random cubic graphs via renewal theory.

Despite this extensive literature, a natural question appears to have remained unanswered: \emph{how many cycles does a uniformly random 2-factor of a random regular graph have?} While the counting argument of \cite{frieze1996} bounds the ratio of the numbers of 2-factors and Hamilton cycles by a polynomial in $n$, it does not determine the typical number of cycles. In this paper, we settle the question exactly for a uniformly random pair $(P,F)$ consisting of a pairing $P$ and a 2-factor $F\subseteq P$, for every $n$ and every $d\ge 3$ simultaneously, and find that the distribution of the cycle count does not depend on $d$.

Throughout, $P$ is a pairing of $n$ vertices of degree $d\ge3$, $G$ is its multigraph (Section~\ref{sec:model}), $F\subseteq P$ is a 2-factor, and $X=X(F)$ is the number of cycles of $F$. Unless $G$ is fixed, the pair $(P,F)$ is uniform over all pairs, and $\E$ and $\Prob$ refer to that probability space.

Our first result is that the distribution of $X$ is explicit and does not depend on $d$ (Theorem~\ref{thm:pgf}): for every $d\ge3$,
\[
\E\,u^{X}\;=\;\frac{\Gamma(n+\tfrac u2)\,\Gamma(\tfrac12)}{\Gamma(n+\tfrac12)\,\Gamma(\tfrac u2)},
\qquad\text{equivalently}\qquad
X\;\stackrel{d}{=}\;1+\sum_{j=1}^{n-1}B_j,
\]
with $B_j\sim\mathrm{Bernoulli}\bigl(\tfrac1{2j+1}\bigr)$ independent. Hence $\E X=\tfrac12\log n+\tfrac\gamma2+\log2+o(1)$, $\operatorname{Var}X=\tfrac12\log n+O(1)$, and $(X-\tfrac12\log n)/\sqrt{\tfrac12\log n}\Rightarrow N(0,1)$ (Corollary~\ref{cor:clt}); this is the analogue for 2-factors of the Feller coupling and of Goncharov's theorem for the cycles of a random permutation \cite{arratia2003,goncharov1944}. If $F$ is required to have girth at least $3$, which is automatic when $G$ is simple, the mean shifts by $-\tfrac34$ and the variance and limiting distribution are unchanged (Theorem~\ref{thm:girth}). In particular (Corollary~\ref{cor:ham}), for every $d\ge3$,
\[
\Prob(X=1)=\tfrac12\sqrt{\pi/n}\,\bigl(1+o(1)\bigr),
\qquad
\Prob(X'=1)=\tfrac{\sqrt\pi\,e^{3/4}}{2}\,n^{-1/2}\bigl(1+o(1)\bigr),
\]
where $X'$ is the cycle count under the girth restriction: one 2-factor in $\Theta(\sqrt n)$ is a Hamilton cycle.

Second, with $Z(u,P)=\sum_F u^{X(F)}$ denoting the generating polynomial of the 2-factors of $P$ by number of cycles, Theorem~\ref{thm:secondmoment} gives $\sum_P Z(u,P)^2$ as an explicit finite double sum for all $n$, $d$, and $u$. Proposition~\ref{prop:second-moment-asymp} (proved in Appendix~\ref{sec:laplace}) shows that $\E Z(u)^2/(\E Z(u))^2\to C(u)$ uniformly on a complex disk about $u=1$, with $C(u)$ given in \eqref{eq:Cu}; $C(1)$ is Robalewska's constant for 2-factors \cite{robalewska1996} and $C(0^+)$ is the Robinson--Wormald constant for Hamilton cycles \cite{robinson_wormald1994}.

Third, the second moment and small subgraph conditioning \cite{janson1995,wormald1999} transfer these results from the uniform pair to a fixed graph. Let $G$ be a uniformly random simple $d$-regular graph and $F$ a uniformly random 2-factor of $G$, and write $\mu_n=\tfrac12\log n+\tfrac\gamma2+\log2-\tfrac34$ and $\sigma_n=\sqrt{\tfrac12\log n}$. Then $\sqrt n\,\Prob(F\text{ is a Hamilton cycle}\mid G)$ converges in distribution to a strictly positive, almost surely finite random variable (Proposition~\ref{prop:typical-ham}), so the expected number of resamplings in the scheme of \cite{frieze1996} is $\Theta_p(\sqrt n)$; $\E[X\mid G]-\mu_n$ converges in distribution to an explicit almost surely finite random variable (Theorem~\ref{thm:typical-mean}); the conditional distribution of $(X-\mu_n)/\sigma_n$ given $G$ converges weakly to $N(0,1)$ in probability (Theorem~\ref{thm:typical-clt}); and $\operatorname{Var}(X\mid G)=\sigma_n^2(1+o_p(1))$ (Proposition~\ref{rem:typical-var}).

Appendix~\ref{sec:numerics} compares these formulas with exhaustive enumeration for small $n$ (for $d=3$, by listing all perfect-matching complements up to $n=80$) and, for $n$ up to $10^6$, with 2-factors produced by the Karp--Sipser algorithm finished by Edmonds' blossom algorithm \cite{karpsipser1981,edmonds1965}.

\section{Pairings, 2-factors, and the probability space}\label{sec:model}
Let $d\ge 3$ be fixed, and let $dn$ be even. In the \emph{pairing model} \cite{bollobas1980,wormald1999}, each vertex $v\in[n]$ owns $d$ points, and a pairing $P$ is a perfect matching of the $dn$ points; contracting the points of each vertex projects $P$ to a $d$-regular multigraph. We write
\[
M(2m)\;=\;\frac{(2m)!}{m!\,2^m}
\]
for the number of perfect matchings on $2m$ points.

A \emph{2-factor of the pairing} $P$ is a subset $F\subseteq P$ of pairs such that every vertex has exactly two of its points covered by $F$. Projected to the multigraph, $F$ is a spanning 2-regular sub-multigraph: a disjoint union of cycles, where a pair joining two points of the same vertex is a \emph{loop} (a cycle of length 1) and two pairs joining the same two vertices form a cycle of length 2. Let $Z_c(P)$ denote the number of 2-factors of $P$ with exactly $c$ cycles, and let $Z^{\ge3}_c(P)$ count only those of girth at least 3 (no loops, no repeated edges); if the multigraph is simple, all its 2-factors have girth at least $3$ automatically. Note that $Z_1 = Z^{\ge 3}_1$ counts Hamilton cycles (for $n\ge3$).

Let $\Omega_{n,d}$ denote the set of pairs $(P,F)$ in which $P$ is a pairing and $F\subseteq P$ is a 2-factor, equipped with the uniform probability $\Prob$. Write $Z(P)=\sum_c Z_c(P)$ for the number of 2-factors of $P$ and
\[
Z(u)\;=\;Z(u,P)\;=\;\sum_c Z_c(P)\,u^c\;=\;\sum_{F} u^{X(F)}
\]
for the generating polynomial of the 2-factors of $P$ by number of cycles, the sum running over the 2-factors $F$ of $P$. Since every 2-factor extends to exactly $M((d-2)n)$ full pairings, the marginal distribution of $P$ under $\Prob$ assigns probability proportional to $Z(P)$ to each pairing, and the cycle count satisfies $\Prob(X=c)=\E Z_c/\E Z$, where here and throughout $\E$ applied to a function of $P$ alone denotes expectation over a uniformly random pairing $P$. Unless we condition on $G$, statements about $X$ refer to this uniform pair. Results for a fixed graph concern the conditional distribution of $X$ given $G$: they are limit theorems, as $n\to\infty$, for the random variables $\Prob(X=1\mid G)$, $\E[X\mid G]$, and $\operatorname{Var}(X\mid G)$, which are functions of $G$, and for the conditional distribution of $X$ given $G$, which converges in probability.

Throughout, we use $L(x) := \log\frac{1}{1-x} = \sum_{k\ge1} x^k/k$.

\section{The exact distribution of the cycle count}\label{sec:exact}

\begin{theorem}\label{thm:moments}
For all $n$, $d\ge3$, and $c\ge1$,
\[
\sum_{P} Z_c(P) \;=\; n!\,[x^n]\,\frac{1}{c!}\left(\frac{L(x)}{2}\right)^{\!c} \cdot \bigl(d(d-1)\bigr)^n \, M\bigl((d-2)n\bigr),
\]
where the sum runs over all $M(dn)$ pairings. Consequently,
\[
\E Z_c
\;=\;
n!\,[x^n]\frac{(L/2)^c}{c!}\,\bigl(d(d-1)\bigr)^n\,
\frac{M((d-2)n)}{M(dn)}.
\]
The same holds for $Z^{\ge3}_c$ with $L(x)/2$ replaced by
\[
L_{\ge3}(x) \;=\; \frac{L(x)}{2}-\frac{x}{2}-\frac{x^2}{4}.
\]
\end{theorem}

\begin{proof}
Both sides count incidences $(P,F)$ where $F$ has $c$ cycles: a partial pairing $F$ covering exactly two points per vertex extends to exactly $M((d-2)n)$ full pairings, so it suffices to count partial pairings with $c$ cycles.

Let $a_k$ be the number of single-cycle partial pairings on a fixed set of $k$ labelled vertices of degree $d$. For $k=1$ (a loop), we choose an unordered pair of the $d$ points: $a_1 = \binom{d}{2} = \tfrac{d(d-1)}{2}$. For $k=2$ (a double edge), we choose two points at each vertex and one of the two crosswise pairings: $a_2 = 2\binom d2^2 = \tfrac{(d(d-1))^2}{2}$. For $k\ge3$, there are $(k-1)!/2$ undirected cyclic orders, and along a fixed traversal, each vertex contributes an ordered choice of (entry point, exit point): $a_k = \tfrac{(k-1)!}{2}\,(d(d-1))^k$. In all cases, $a_k\,\frac{x^k}{k!} = \tfrac{(d(d-1)x)^k}{2k}$, so the exponential generating function of single cycles is $\tfrac12 L\bigl(d(d-1)x\bigr)$. By the exponential formula, the number of partial pairings with exactly $c$ cycles on $[n]$ is
\[
n!\,[x^n]\,\frac{1}{c!}\Bigl(\tfrac12 L\bigl(d(d-1)x\bigr)\Bigr)^{c}
= \bigl(d(d-1)\bigr)^n\, n!\,[x^n]\,\frac{(L(x)/2)^c}{c!}.
\]
Excluding $k\in\{1,2\}$ removes the terms $x/2$ and $x^2/4$ from $L/2$, yielding the girth-$\ge3$ statement.
\end{proof}

\begin{corollary}\label{cor:small-c}
For all $n\ge 1$ and $d\ge 3$, we have $\E Z_1 = \frac{(n-1)!}{2}\,(d(d-1))^n \frac{M((d-2)n)}{M(dn)}$, which recovers the expected number of Hamilton cycles in \cite{frieze1996}. Moreover, $n!\,[x^n]\frac{(L(x)/2)^2}{2} = \frac{(n-1)!}{4}H_{n-1}$, where $H_m=\sum_{k=1}^{m}1/k$ is the $m$-th harmonic number, and therefore $\E Z_2 / \E Z_1 = \tfrac{1}{2}H_{n-1}$.
\end{corollary}

\begin{proof}
We have $[x^n]L = 1/n$ and $[x^n]L^2 = \sum_{k=1}^{n-1}\frac{1}{k(n-k)} = \frac{2}{n}H_{n-1}$.
\end{proof}

\begin{theorem}\label{thm:pgf}
For a uniform pair $(P,F)$ and every $d \ge 3$,
\[
\E\,u^{X} \;=\; \frac{[y^n]\,(1-y)^{-u/2}}{[y^n]\,(1-y)^{-1/2}}
\;=\; \frac{\Gamma\!\left(n+\tfrac u2\right)\Gamma\!\left(\tfrac 12\right)}{\Gamma\!\left(n+\tfrac 12\right)\Gamma\!\left(\tfrac u2\right)}
\;=\; u\prod_{j=1}^{n-1}\frac{u+2j}{1+2j}.
\]
In particular, the distribution of $X=X(F)$ does not depend on $d$, and
\[
X \;\stackrel{d}{=}\; 1+\sum_{j=1}^{n-1} B_j,\qquad B_j\sim\mathrm{Bernoulli}\Bigl(\frac{1}{2j+1}\Bigr)\ \text{independent.}
\]
\end{theorem}

\begin{proof}
Multiplying the identity of Theorem~\ref{thm:moments} by $u^c$ and summing over $c\ge 1$, we obtain from $\sum_{c\ge 0} \frac{(uL(x)/2)^c}{c!} = e^{uL(x)/2} = (1-x)^{-u/2}$ that
\[
\sum_{c\ge 1} u^c \sum_P Z_c(P) = n!\,[x^n](1-x)^{-u/2}\,(d(d-1))^n M\bigl((d-2)n\bigr).
\]
Since $\E u^X$ is the ratio of this sum to its value at $u=1$, all factors depending on $d$ cancel, which gives the first equality. The second equality follows from the binomial expansion
\[
[y^n](1-y)^{-u/2} = \binom{n+u/2-1}{n} = \frac{\Gamma(n+u/2)}{\Gamma(u/2)\,n!},
\]
and the product form follows from $\Gamma(n+u/2)/\Gamma(u/2) = \prod_{j=0}^{n-1}(j+u/2) = 2^{-n}u\prod_{j=1}^{n-1}(u+2j)$ upon dividing by its value at $u=1$. Finally, the prefactor $u$ corresponds to the constant $1$, and each factor $\frac{u+2j}{1+2j}$ is the probability generating function of an independent $\mathrm{Bernoulli}\bigl(\frac1{2j+1}\bigr)$ random variable.
\end{proof}

\begin{corollary}\label{cor:clt}
For a uniform pair $(P,F)$ and every $d\ge 3$, we have $\E X = H_{2n-1}-\tfrac12 H_{n-1} = \tfrac12\log n + \tfrac{\gamma}{2} + \log 2 + o(1)$,
$\operatorname{Var}X = \tfrac12\log n + O(1)$, and
\[
\frac{X-\tfrac12\log n}{\sqrt{\tfrac12 \log n}}\;\Longrightarrow\; N(0,1).
\]
\end{corollary}

\begin{proof}
The mean and variance follow immediately from the decomposition of $X$ into independent Bernoulli variables in Theorem~\ref{thm:pgf}, and the central limit theorem follows from Lindeberg's condition, just as in Goncharov's theorem for the number of cycles of a random permutation \cite{arratia2003,goncharov1944}.
\end{proof}

\begin{corollary}\label{cor:quasipower-input}
Fix $d\ge 3$ and $\varepsilon\in(0,1)$. Let
\begin{align*}
m_n
&\;:=\;
\E X
\;=\;
H_{2n-1}-\tfrac12 H_{n-1}
\;=\;
\tfrac12\log n+\tfrac\gamma2+\log2+o(1),\\
v_n
&\;:=\;
\operatorname{Var}X
\;=\;
\sum_{j=1}^{n-1}\frac{2j}{(2j+1)^2}
\;=\;
\tfrac12\log n+O(1)
\end{align*}
be the exact pairing mean and variance from Corollary~\ref{cor:clt}, so that $m_n=\mu_n^{\mathrm{pair}}+o(1)$ and $v_n=\sigma_n^2+O(1)$, where $\mu_n^{\mathrm{pair}}=\tfrac12\log n+\tfrac\gamma2+\log2$ and $\sigma_n=\sqrt{\tfrac12\log n}$. Then, uniformly for $s\in\mathbb C$ with $|s|\le\varepsilon$, we have
\begin{equation}\label{eq:EZ-complex}
\frac{\E Z(e^{s})}{\E Z(1)}
\;=\;
\exp\Bigl(m_n s+\tfrac12 v_n s^2+O(|s|^3\log n)\Bigr).
\end{equation}
In terms of the additive parameter $t=u-1$, all coefficients beyond the linear term are $O(1)$:
\[
\frac{\E Z(1+t)}{\E Z(1)}
\;=\;
\exp\Bigl(m_n t-\frac{\pi^2}{16}t^2+O(|t|^3)\Bigr).
\]
In particular, at $u=e^{it/\sigma_n}$ the characteristic function of $X$ for a uniform pair satisfies the Gaussian expansion of Hwang's quasi-power theorem \cite{hwang1998}.
\end{corollary}

\begin{proof}
Since every 2-factor extends to $M((d-2)n)$ pairings, $\E Z(u)/\E Z(1)=\E u^X=\Gamma(n+u/2)\Gamma(1/2)/(\Gamma(n+1/2)\Gamma(u/2))$ by Theorem~\ref{thm:pgf}. Equivalently, $X\stackrel{d}{=}1+\sum_{j=1}^{n-1}B_j$ with independent $B_j\sim\mathrm{Bernoulli}(1/(2j+1))$, so the cumulant generating function in $s=\log u$ is
\[
\log\E e^{sX}
\;=\;
s+\sum_{j=1}^{n-1}\log\bigl(1-p_j+p_j e^{s}\bigr),\qquad p_j=\frac1{2j+1}.
\]
For $|s|\le\varepsilon$, we have the uniform expansion $\log(1-p_j+p_j e^{s})=p_j s+\tfrac12 p_j(1-p_j)s^2+O(p_j|s|^3)$. Summing over $1\le j\le n-1$ gives \eqref{eq:EZ-complex}, since $\E X=1+\sum_j p_j=m_n$, $\operatorname{Var}X=\sum_jp_j(1-p_j)=v_n$, and $\sum_{j<n}p_j=O(\log n)$ (the third cumulant of $X$ is $\sum_j p_j(1-p_j)(1-2p_j)\sim\tfrac12\log n$). In terms of $t=u-1$, differentiating $\log\E Z(u)$ twice at $u=1$ gives $\partial_{uu}\log\E Z(u)|_{u=1}=\tfrac14\psi'(n+\tfrac12)-\tfrac14\psi'(\tfrac12)=-\pi^2/8+O(n^{-1})$ because $\psi'(1/2)=\pi^2/2$, which yields the coefficient $-\pi^2/16$ of $t^2$.
\end{proof}

\begin{remark}
The independence of $d$ in Theorem~\ref{thm:pgf} has a direct combinatorial explanation: once we specify which two of the $d$ points at each vertex belong to the 2-factor, the remaining $(d-2)n$ points can be paired in $M((d-2)n)$ ways regardless of how the $2n$ chosen points are paired. Consequently, under the uniform distribution on pairs $(P,F)$, the induced factor $F$ is a uniformly random pairing of $2n$ points with two points at each vertex, and its cycle structure is that of a random 2-regular multigraph on $[n]$.
\end{remark}

\section{Girth at least three, and the probability of a Hamilton cycle}\label{sec:girth}

A simple host graph has no loops or digons, so all of its 2-factors have girth at least $3$. The girth-$\ge3$ distribution below is therefore the one relevant to simple $G(n,d)$; Lemma~\ref{lem:simple-transfer} makes this precise.

\begin{theorem}\label{thm:girth}
Let $X'$ be the number of cycles of $F$ when $(P,F)$ is uniform over the pairs in which $F$ has girth at least $3$. Then, uniformly for $u$ in compact subsets of $\mathbb C\setminus\{0,-2,-4,\dots\}$,
\[
\E\,u^{X'} \;=\; \frac{[y^n]\,e^{-u(y/2+y^2/4)}(1-y)^{-u/2}}{[y^n]\,e^{-(y/2+y^2/4)}(1-y)^{-1/2}}
\;=\; e^{-\frac{3(u-1)}{4}}\, n^{\frac{u-1}{2}}\,\frac{\Gamma(1/2)}{\Gamma(u/2)}\,\bigl(1+o(1)\bigr).
\]
Consequently, $\E X' = \tfrac12\log n + \tfrac{\gamma}{2}+\log 2 - \tfrac34 + o(1)$, $\operatorname{Var}X' = \tfrac12\log n + O(1)$, and $X'$ satisfies the same central limit theorem as $X$.
\end{theorem}

\begin{proof}
The first equality follows from Theorem~\ref{thm:moments} summed against $u^c$, with $e^{uL_{\ge3}(x)} = e^{-u(x/2+x^2/4)}(1-x)^{-u/2}$. Since $e^{-u(y/2+y^2/4)}$ is entire, standard singularity analysis \cite[Ch.~VI]{flajolet2009} yields
$[y^n] e^{-u(y/2+y^2/4)}(1-y)^{-u/2} = e^{-3u/4}\, n^{u/2-1}/\Gamma(u/2)\,(1+O(1/n))$, uniformly on compact subsets of $\mathbb C\setminus\{0,-2,-4,\dots\}$, and in particular on a complex neighbourhood of $u=1$ (multiplying the expansions of Corollary~\ref{cor:quasipower-input} by $e^{-3(u-1)/4}$). The moment asymptotics and the CLT then follow from the quasi-power theorem of Hwang \cite{hwang1998}.
\end{proof}

\begin{corollary}\label{cor:ham}
For a uniform pair $(P,F)$,
\begin{align*}
\Prob(X=1)
&\;=\;
\frac{2^{2n-1}\,n!\,(n-1)!}{(2n)!}
\;=\;
\frac12\sqrt{\frac{\pi}{n}}\,\bigl(1+o(1)\bigr),\\
\Prob(X'=1)
&\;=\;
\frac{\sqrt{\pi}\,e^{3/4}}{2}\,n^{-1/2}\,\bigl(1+o(1)\bigr),
\end{align*}
for every $d\ge3$. More generally, for each fixed $c$,
$\Prob(X'=c) = \frac{\sqrt{\pi}\,e^{3/4}}{2\sqrt n}\,\frac{\lambda^{c-1}}{(c-1)!}(1+o(1))$ with $\lambda=\tfrac12\log n$.
\end{corollary}

\begin{proof}
By Theorem~\ref{thm:pgf}, we have $\Prob(X=1) = \prod_{j=1}^{n-1}\frac{2j}{2j+1} = \frac{2^{2n-1}n!(n-1)!}{(2n)!}$, and Stirling's formula gives the asymptotic estimate. For $X'$, we have $\Prob(X'=c) = N_c/\sum_c N_c$, where $N_c = n![x^n]L_{\ge3}(x)^c/c!$; singularity analysis yields $N_c \sim \frac{(n-1)!}{2}\frac{(\frac12\log n)^{c-1}}{(c-1)!}$ for fixed $c$ and $\sum_c N_c = n![x^n]e^{L_{\ge3}(x)} \sim n!\,e^{-3/4}/\sqrt{\pi n}$.
\end{proof}

\begin{remark}\label{rem:fjmrw}
Corollary~\ref{cor:ham} states that, on average, one in $\Theta(\sqrt n)$ uniform 2-factors is Hamiltonian, with the explicit constant $\frac{2}{\sqrt\pi}$ in the ratio $\sum_c\E Z_c/\E Z_1 = \frac{2}{\sqrt{\pi}}\sqrt n\,(1+o(1))$. Combining the bounds $Z \le n\,\E Z$ (by Markov's inequality) and $Z_1 \ge \E Z_1/n$ (via small subgraph conditioning \cite{janson1995,wormald1999}), as in \cite{frieze1996}, gives $\Prob(F\text{ Hamiltonian}\mid G)=\Omega(n^{-5/2})$ with high probability. Proposition~\ref{prop:typical-ham} sharpens this bound to $\Prob(F\text{ Hamiltonian}\mid G) = (\Lambda/\sqrt n)(1+o_p(1)) = \Theta_p(n^{-1/2})$, showing that for a uniformly random simple $G$ the expected number of uniform 2-factor resamplings needed to obtain a Hamilton cycle is $1/\Prob(F\text{ Hamiltonian}\mid G)=\sqrt n/\Lambda\,(1+o_p(1))=\Theta_p(\sqrt n)$.
\end{remark}

\section{Short cycles of $G$ and transfer to simple graphs}\label{sec:ssc}

Let $X_i = X_i(G)$ denote the number of $i$-cycles of the host multigraph $G$ (not to be confused with $X=X(F)$, the cycle count of the 2-factor). In the pairing model, the variables $X_i$ converge in distribution to independent Poisson variables with means $\lambda_i = (d-1)^i/(2i)$ \cite{bollobas1980,wormald1999}. Small subgraph conditioning \cite{janson1995} states that for a count $Y$ of spanning structures satisfying $\E[Y X_i]/\E Y \to \lambda_i(1+\delta_i)$ together with $\E Y^2/(\E Y)^2\to\exp(\sum_{i\ge1}\lambda_i\delta_i^2)$, one has
\[
\frac{Y}{\E Y} \;\longrightarrow\; W_Y \;=\; \prod_{i\ge1}(1+\delta_i)^{X_i}e^{-\lambda_i\delta_i}
\]
in the sense that the ratio $(Y/\E Y)/W_Y$ converges to $1$ in probability. The two basic counts we need are $Z_1$ (Hamilton cycles; Robinson--Wormald \cite{robinson_wormald1994}, cf.\ \cite{frieze1996,janson1995}) and $Z=\sum_c Z_c$ (2-factors; Robalewska \cite{robalewska1996}). The corresponding parameters $\delta_i$, as well as their extension to the generating polynomial $Z(u)=\sum_F u^{X(F)}$, admit a unified derivation. (For real $u>0$, the ratio $\E[Z(u)f]/\E Z(u)$ is the expectation of $f(P)$ when $(P,F)\in\Omega_{n,d}$ is drawn with probability proportional to $u^{X(F)}$.)

\begin{lemma}\label{lem:delta}
For a uniform pair $(P,F)\in\Omega_{n,d}$, conditioning on $F$ leaves $P=F\cup R$, where $R$ is a uniformly random pairing of the remaining $(d-2)n$ points, and the expected number of $i$-cycles of the resulting multigraph satisfies
\[
\E\, X_i(F\cup R) \;\longrightarrow\; \frac{(d-1)^i+(-1)^i}{2i},
\qquad\text{hence}\qquad
\delta_i^{2\mathrm F} = \frac{(-1)^i}{(d-1)^i}.
\]
If $F$ is constrained to be a single Hamilton cycle, the all-$F$ contribution (an $i$-cycle that forms a component of $F$) vanishes, yielding
\[
\delta_i^{\mathrm H} = \frac{(-1)^i-1}{(d-1)^i},
\]
which vanishes for even $i$, recovering the well-known fact that only \emph{odd} short cycles affect Hamilton cycle counts \cite{robinson_wormald1994}.
More generally, for fixed $i\ge 1$,
\begin{equation}\label{eq:delta-u-lem}
\frac{\E\bigl[Z(u)\,X_i\bigr]}{\E Z(u)}\;\longrightarrow\;
\lambda_i\bigl(1+\delta_i(u)\bigr),\qquad
\delta_i(u)=\frac{u-1+(-1)^i}{(d-1)^i},
\end{equation}
with $\lambda_i=(d-1)^i/(2i)$. In particular, $\delta_i(1)=\delta_i^{2\mathrm F}$ and $\delta_i(u)\to\delta_i^{\mathrm H}$ as $u\to 0^+$.
\end{lemma}

\begin{proof}
For a uniform pair, conditioning on $F$ leaves $R$ as a uniform pairing of the remaining $(d-2)n$ points, so $G=F\cup R$. Fix $i\ge 1$ and consider the expected number of directed closed walks of length $i$ that project to undirected $i$-cycles of $G$, then divide by $2i$.

Classify each edge of a prospective $i$-cycle as an $F$-edge or an $R$-edge. At a vertex whose two incident cycle edges are both of $F$-type, the ordered choice of (entry, exit) points contributes a factor of $2$; when both are of $R$-type, the factor is $(d-2)(d-3)$; and in the mixed case, it is $2(d-2)$. Each $F$-edge has pairing probability $\sim 1/(2n)$ among the $2n$ points assigned to $F$, and each $R$-edge has pairing probability $\sim 1/((d-2)n)$ among the remaining points. Recording the sequence of edge types around the cycle via the $2\times 2$ transition matrix
\[
T=\begin{pmatrix}1 & \sqrt{2(d-2)}\\[2pt] \sqrt{2(d-2)} & d-3\end{pmatrix},
\]
we find that the eigenvalues of $T$ are $d-1$ and $-1$, so
\[
\E X_i(F\cup R)\;\longrightarrow\;\frac{\operatorname{tr}(T^i)}{2i}=\frac{(d-1)^i+(-1)^i}{2i}.
\]
Dividing by $\lambda_i=(d-1)^i/(2i)$ yields $1+\delta_i^{2\mathrm F}$.

When $F$ is conditioned to be a Hamilton cycle, the configuration in which the entire $i$-cycle lies in $F$ is precluded (a Hamilton cycle consists of a single component of length $n$). Removing that contribution subtracts $1/(2i)$ from the limit, yielding
\[
\frac{(d-1)^i+(-1)^i-1}{2i}=\lambda_i\bigl(1+\delta_i^{\mathrm H}\bigr),
\]
in agreement with \cite[Lemma~2]{frieze1996}.

When each 2-factor $F$ is weighted by $u^{X(F)}$, any $i$-cycle that uses at least one edge of $R$ intersects $F$ in proper paths, constraining only $O(1)$ points and leaving the expected weight of the completion of $F$ unchanged up to a factor of $1+O(i/n)$. By contrast, an $i$-cycle contained entirely in $F$ forms an additional component of $F$ and is multiplied by the weight $u$. Hence $\E[Z(u)X_i]/\E Z(u)$ converges to
\[
\frac{(d-1)^i+(-1)^i-1+u}{2i}=\lambda_i\bigl(1+\delta_i(u)\bigr),
\]
which proves \eqref{eq:delta-u-lem}.
\end{proof}

\begin{lemma}\label{lem:factorial}
Fix $d\ge 3$, let $i_1,\dots,i_m\ge 1$ and $r_1,\dots,r_m\ge 1$ be integers, and write $(x)_r=x(x-1)\cdots(x-r+1)$. Then
\[
\frac{\E\Bigl[Z(u)\prod_{j=1}^{m}(X_{i_j})_{r_j}\Bigr]}{\E Z(u)}
\;\longrightarrow\;
\prod_{j=1}^{m}\lambda_{i_j}^{r_j}\bigl(1+\delta_{i_j}(u)\bigr)^{r_j}
\]
uniformly for $u$ in a sufficiently small complex neighbourhood of $1$ (and in particular along every sequence $u_n\to1$). In particular, for each fixed $i$ and $r\ge 1$,
\[
\frac{\E\bigl[Z(u)\,(X_i)_r\bigr]}{\E Z(u)}\;\longrightarrow\;\lambda_i^r\bigl(1+\delta_i(u)\bigr)^r.
\]
\end{lemma}

\begin{proof}
As in Lemma~\ref{lem:delta}, $G=F\cup R$. The falling factorial $(X_i)_r$ enumerates ordered $r$-tuples of distinct $i$-cycles of $G$. For fixed $i$ and $r$, mark $r$ vertex-disjoint directed $i$-cycles and determine their contribution by the same transition-matrix computation: each marked cycle contributes a factor of $\lambda_i(1+\delta_i(u))$ in the limit, because (i) a cycle using at least one edge of $R$ constrains only proper paths of $F$, and the underlying Gamma ratios $\Gamma(n-k+u/2)\,n!/\bigl(\Gamma(n+u/2)\,(n-k)!\bigr)$ converge to $1$ uniformly on a complex neighbourhood of $u=1$ up to an error of $O(i/n)$, and (ii) a short cycle contained entirely in $F$ is weighted by $u$, as in the case $r=1$ of Lemma~\ref{lem:delta}. Configurations in which two marked cycles share a vertex contribute $O(n^{-1})$ after normalisation by $\E Z(u)$: the union of two distinct $i$-cycles with a common vertex has $v\le2i-1$ vertices and $e\ge v+1$ edges, and every such configuration carries weight $O(n^{v-e})$ (a factor $O(n)$ per vertex label, $O(1)$ per point choice, and $O(1/n)$ per pair), whence
\begin{align*}
&\frac{\E\bigl[Z(u)\,\#\{\text{ordered pairs of distinct $i$-cycles sharing a vertex}\}\bigr]}{\E Z(u)}\\
&\qquad=\;
\sum_{v\le 2i-1}\sum_{e\ge v+1}O\bigl(n^{v-e}\bigr)
\;=\;
O(n^{-1}).
\end{align*}
Thus, the contributions of the $r$ marked cycles multiply, yielding $\lambda_i^r(1+\delta_i(u))^r$.

For several lengths $i_1,\dots,i_m$, mark $r_j$ cycles of each length $i_j$. As long as the total number of marked vertices remains $O(1)$, the same computation applies simultaneously: cross-terms in which distinct marked cycles share vertices are negligible, and the exponential formula factorises the contributions, yielding the product on the right-hand side. This extends Robalewska's first-moment asymptotic $\E[Z\,X_i]/\E Z\to\lambda_i(1+\delta_i)$ \cite[Lemma~1, Corollary~1]{robalewska1996} from a single short cycle to falling factorials and general $u$.
\end{proof}

\begin{lemma}\label{lem:simple-transfer}
Fix $d\ge 3$. Conditioning the pairing model on the underlying multigraph being simple transfers every event of pairing-probability $o(1)$ to an event of simple-graph probability $o(1)$. When the underlying multigraph is simple, every 2-factor has girth at least $3$, so $X=X'$. Consequently, conditioned on simplicity, the distribution of $X$ for a uniform pair is asymptotically the girth-$\ge3$ distribution of Theorem~\ref{thm:girth}:
\[
\E\bigl[u^{X}\bigm|\text{simple}\bigr]
\;=\;
e^{-3(u-1)/4}\,\E u^{X}\,\bigl(1+o(1)\bigr)
\]
uniformly for $u$ in a neighbourhood of $1$, where $\E u^{X}$ is the generating function of Theorem~\ref{thm:pgf}. In particular, the mean, variance, and central limit theorem of a uniformly random 2-factor of a uniformly random \emph{simple} $d$-regular graph are those of Theorem~\ref{thm:girth}, not those of Corollary~\ref{cor:clt}.
\end{lemma}

\begin{proof}
By the classical enumeration of regular graphs \cite{bollobas1980,wormald1999}, $\Prob(\text{simple})\to c_d\in(0,1)$, a constant depending only on $d$. Hence, for any event $A_n$ with $\Prob(A_n)\to0$ under the unconditioned pairing model,
\[
\Prob(A_n\mid\text{simple})
\;=\;
\frac{\Prob(A_n\cap\text{simple})}{\Prob(\text{simple})}
\;\le\;
\frac{\Prob(A_n)}{c_d-o(1)}
\;\longrightarrow\; 0.
\]
This suffices for all subsequent conditional $o_p(1)$ statements.

For the distribution of $X$, a simple host graph has no loops or digons, so every 2-factor has girth at least $3$ and $X=X'$ (which gives the generating-function asymptotic directly for complex $u$ by Theorem~\ref{thm:girth}). From the short-cycle perspective for real $u>0$, simplicity of the pairing is the event $\{X_1(G)=X_2(G)=0\}$. By Lemma~\ref{lem:delta}, $\E[Z(u)X_i]/\E Z(u)\to\lambda_i(1+\delta_i(u))$ for $i=1,2$, with $\lambda_i=(d-1)^i/(2i)$ and $\delta_i(u)=(u-1+(-1)^i)/(d-1)^i$; the corresponding factorial moments (Lemma~\ref{lem:factorial}) imply that under the $Z(u)$-weighted pairing distribution, $(X_1,X_2)$ converge in distribution to independent Poisson variables with means $\lambda_i(1+\delta_i(u))$. Therefore,
\[
\frac{\E\bigl[Z(u)\,1_{\mathrm{simple}}\bigr]}{\E Z(u)}
\;\longrightarrow\;
\exp\bigl(-\lambda_1(1+\delta_1(u))-\lambda_2(1+\delta_2(u))\bigr),
\]
and the conditional generating function is
\begin{align*}
\E\bigl[u^{X}\bigm|\text{simple}\bigr]
&\;=\;
\frac{\E\bigl[Z(u)\,1_{\mathrm{simple}}\bigr]}{\E\bigl[Z(1)\,1_{\mathrm{simple}}\bigr]}\\
&\;=\;
\frac{\E Z(u)}{\E Z(1)}\,
\exp\Bigl(-\sum_{i=1,2}\lambda_i\bigl(\delta_i(u)-\delta_i(1)\bigr)\Bigr)\,(1+o(1)).
\end{align*}
The exponent is $\lambda_1(\delta_1(u)-\delta_1(1))+\lambda_2(\delta_2(u)-\delta_2(1))=(u-1)/2+(u-1)/4=3(u-1)/4$, matching the multiplicative factor $e^{-3(u-1)/4}$ of Theorem~\ref{thm:girth}.
\end{proof}

\section{The second moment of $Z(u)$}\label{sec:variance}

For 2-factors --- unlike Hamilton cycles --- the second moment admits an \emph{exact} finite-$n$ formula, because the overlap of two 2-factors requires no connectivity accounting. Write $(q)_m = q(q+1)\cdots(q+m-1)$, and recall $Z(u)=Z(u,P)=\sum_F u^{X(F)}$ from Section~\ref{sec:model}.

\begin{theorem}\label{thm:secondmoment}
For all $n$, $d\ge3$, and $u$, with $a = d(d-1)$, $b=(d-2)(d-3)$, and $c_e = d(d-1)(d-2)$,
\begin{align*}
\sum_P Z(u,P)^2
&\;=\;
n!\sum_{n_0,\,p\ \ge0}
\frac{\bigl((\tfrac u2)_{n_0+p}\bigr)^{2}\,(ab)^{n_0}\,4^{p}}{n_0!\;p!}\\
&\qquad\times
[x^{\,n-n_0}]\Bigl(e^{\frac{u^2}{2}L(ax)}\,\hat\pi(x)^{p}\Bigr)\,
M\bigl((d-4)n+2s\bigr),
\end{align*}
where $\hat\pi(x) = \dfrac{c_e^{2}\,x^{2}}{2\,(1-ax)}$ and $s = n-n_0-p$.
\end{theorem}

\begin{proof}
We count triples $(P,F_1,F_2)$ weighted by $u^{X(F_1)+X(F_2)}$. Colour the pairs of $F_1\setminus F_2$ red, the pairs of $F_2\setminus F_1$ blue, and the shared pairs $F_1\cap F_2$ purple. The purple pairs form vertex-disjoint \emph{purple} paths and cycles (corresponding to Robalewska's classification of $R_1\cap R_2$ into paths, cycles, and crosses \cite[\S2]{robalewska1996}, where our free vertices are her crosses). Since every point of the pairing has a unique partner, sharing a point is equivalent to sharing a pair, so the red ($F_1$-only) and blue ($F_2$-only) edges use disjoint point sets and can never coincide. Each vertex is of one of four types, with the indicated number of ways to assign roles to its points: \emph{free} (two red points and two blue points, all distinct: $d(d-1)(d-2)(d-3)=ab$), a \emph{purple-path endpoint} (one purple, one red, and one blue point: $c_e$), or a \emph{purple-path interior} or \emph{purple-cycle} vertex (two purple points: $a$). Figure~\ref{fig:purple} displays the local configurations in the pairing: a purple path versus a purple cycle, followed by the endpoint split and the free cross.

\begin{figure}[H]
\centering
\begin{tikzpicture}[
  vtx/.style={circle, draw=black!75, fill=white, inner sep=1.45pt, minimum size=5.2pt},
  F1/.style={
    draw=none,
    decoration={markings,
      mark=between positions 0.06 and 0.94 step 2.15pt with {
        \pgftransformresetnontranslations
        \draw[red!75!black, line width=0.5pt, cap=round]
          (-0.95pt,0.95pt) -- (0.95pt,-0.95pt);
      }},
    postaction={decorate},
  },
  F2/.style={
    draw=none,
    decoration={markings,
      mark=between positions 0.06 and 0.94 step 2.15pt with {
        \pgftransformresetnontranslations
        \draw[blue!70!black, line width=0.5pt, cap=round]
          (-0.95pt,-0.95pt) -- (0.95pt,0.95pt);
      }},
    postaction={decorate},
  },
  purple/.style={
    draw=none,
    decoration={markings,
      mark=between positions 0.06 and 0.94 step 2.15pt with {
        \pgftransformresetnontranslations
        \draw[violet!80!black, line width=0.48pt, cap=round]
          (-0.95pt,0.95pt) -- (0.95pt,-0.95pt);
        \draw[violet!80!black, line width=0.48pt, cap=round]
          (-0.95pt,-0.95pt) -- (0.95pt,0.95pt);
      }},
    postaction={decorate},
  },
  lab/.style={font=\scriptsize, text=black!70},
  tag/.style={font=\scriptsize\bfseries, text=black!80, anchor=south west},
  F1a/.style={
    draw=none,
    decoration={markings,
      mark=between positions 0.08 and 0.92 step 3.1pt with {
        \draw[red!75!black, line width=0.55pt, cap=round]
          (0,-1.2pt) -- (0,1.2pt);
      }},
    postaction={decorate},
  },
  F2a/.style={
    draw=none,
    decoration={markings,
      mark=between positions 0.08 and 0.92 step 3.1pt with {
        \draw[blue!70!black, line width=0.55pt, cap=round]
          (0,-1.2pt) -- (0,1.2pt);
      }},
    postaction={decorate},
  },
  purplea/.style={
    draw=none,
    decoration={markings,
      mark=between positions 0.06 and 0.94 step 4.2pt with {
        \pgftransformresetnontranslations
        \draw[violet!80!black, line width=0.48pt, cap=round]
          (-0.72pt,0.72pt) -- (0.72pt,-0.72pt);
        \draw[violet!80!black, line width=0.48pt, cap=round]
          (-0.72pt,-0.72pt) -- (0.72pt,0.72pt);
      }},
    postaction={decorate},
  },
]
\begin{scope}[shift={(0,0)}]
  \coordinate (L) at (0.05,0.85);
  \coordinate (I) at (1.35,0.85);
  \coordinate (R) at (2.65,0.85);
  \coordinate (T) at (1.35,1.78);
  \coordinate (Bo) at (1.35,-0.08);
  \path[purplea] (L) -- (I) -- (R);
  \path[F1a] (L) -- (T) -- (R);
  \path[F2a] (L) -- (Bo) -- (R);
  \node[vtx] at (L) {};
  \node[vtx] at (I) {};
  \node[vtx] at (R) {};
  \node[vtx] at (T) {};
  \node[vtx] at (Bo) {};
  \node[lab, text=red!75!black] at (2.15,1.78) {$F_1$};
  \node[lab, text=blue!70!black] at (2.15,-0.08) {$F_2$};
  \node[lab] at (0.05,0.52) {end};
  \node[lab] at (1.35,0.52) {interior};
  \node[lab] at (2.65,0.52) {end};
\end{scope}
\begin{scope}[shift={(3.55,0.15)}]
  \node[vtx] (P1) at (0.28,1.55) {};
  \node[vtx] (P2) at (1.62,1.55) {};
  \node[vtx] (P3) at (1.62,0.35) {};
  \node[vtx] (P4) at (0.28,0.35) {};
  \path[purple] (P1) -- (P2) -- (P3) -- (P4) -- cycle;
  \node[lab] at (0.95,-0.18) {$\textcolor{red!75!black}{F_1}\mathrel{\textcolor{black!70}{=}}\textcolor{blue!70!black}{F_2}$};
\end{scope}
\begin{scope}[shift={(6.15,0.05)}]
  \node[vtx] (Q) at (0.85,1.55) {};
  \node[vtx] (E) at (0.85,0.72) {};
  \node[vtx] (Rr) at (0.12,-0.05) {};
  \node[vtx] (Bb) at (1.58,-0.05) {};
  \path[purple] (Q) -- (E);
  \path[F1] (E) -- (Rr);
  \path[F2] (E) -- (Bb);
  \node[lab] at (0.85,-0.42) {colours split};
\end{scope}
\begin{scope}[shift={(8.55,0.05)}]
  \node[vtx] (W) at (0.05,0.72) {};
  \node[vtx] (Fv) at (0.8,0.72) {};
  \node[vtx] (Ea) at (1.55,0.72) {};
  \node[vtx] (N) at (0.8,1.55) {};
  \node[vtx] (S) at (0.8,-0.11) {};
  \path[F1] (W) -- (Fv) -- (Ea);
  \path[F2] (N) -- (Fv) -- (S);
  \node[lab] at (0.80,-0.48) {no shared edge};
\end{scope}
\node[tag] at (0.00,-1.02) {a.};
\node[tag] at (3.55,-1.02) {b.};
\node[tag] at (6.15,-1.02) {c.};
\node[tag] at (8.55,-1.02) {d.};
\begin{scope}[shift={(0,-1.48)}]
  \path[purple] (0,0) -- (0.72,0);
  \node[lab, anchor=west] at (0.80,0) {shared pair ($F_1\cap F_2$)};
  \path[F1] (4.25,0) -- (4.97,0);
  \node[lab, anchor=west] at (5.05,0) {$F_1$ only};
  \path[F2] (6.75,0) -- (7.47,0);
  \node[lab, anchor=west] at (7.55,0) {$F_2$ only};
\end{scope}
\end{tikzpicture}
\caption{Overlap of two 2-factors $F_1$ (red) and $F_2$ (blue) of the same pairing $P$. Shared pairs $F_1\cap F_2$ are drawn in purple and form vertex-disjoint paths and cycles. (a)~A cycle of $F_1$ and a cycle of $F_2$ may share a purple path without coinciding. (b)~If every pair of a cycle lies in both factors, that cycle is purple (weight $u^2$). (c)~At a purple-path endpoint, one purple pair meets one red ($F_1$-only) pair and one blue ($F_2$-only) pair. (d)~A free vertex (Robalewska's ``cross'', for $d\ge4$) lies on two red pairs and two blue pairs.}
\label{fig:purple}
\end{figure}
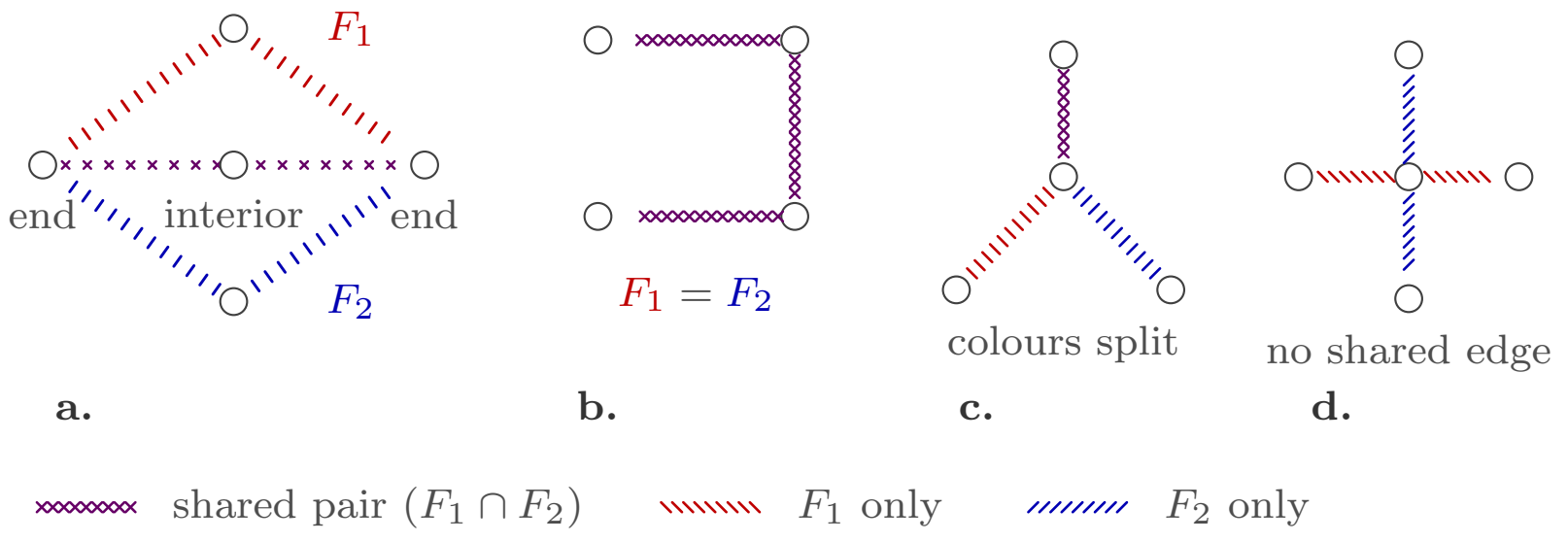

Purple cycles are cycles shared by both factors and carry weight $u^2$ each; by the computation in the proof of Theorem~\ref{thm:moments}, their exponential generating function is $\exp\bigl(\tfrac{u^2}{2}L(ax)\bigr)$. A purple path with $j\ge1$ edges spans $j+1$ vertices (ordered up to reversal) and admits $c_e^2 a^{j-1}$ choices of points, yielding the generating function $\hat\pi(x)$ per path.

Given the purple structure, contract each purple path to a single vertex equipped with two unpaired red ($F_1$) points and two unpaired blue ($F_2$) points (one at each endpoint, with a factor of $2$ for the traversal orientation in each colour). The red completion is then a $u$-weighted 2-regular pairing on the $n_0+p$ vertices, each cycle of which expands to a single cycle of $F_1$; the same holds for the blue completion of $F_2$, \emph{independently}, since the two colours draw from disjoint sets of points. For $n_0$ vertices of weight $\alpha$ and $p$ vertices of weight $\beta$, the exponential formula yields
\[
\sum_{\text{2-regular pairings }\sigma} u^{X(\sigma)}
= n_0!\,p!\,[y^{n_0}\zeta^{p}]\,(1-\alpha y-\beta\zeta)^{-u/2}
= \Bigl(\tfrac u2\Bigr)_{n_0+p}\alpha^{n_0}\beta^{p},
\]
with $(\alpha,\beta)=(a,2)$ for the red completion and $(b,2)$ for the blue completion (the vertex factor $ab$ of a free vertex splitting as $a\cdot b$ between the two colours). Finally, $F_1\cup F_2$ covers $4n-2s$ points, where $s = w+e = n-n_0-p$ is the number of purple pairs ($w$ purple-cycle vertices and $e$ purple-path edges), so it extends to $M((d-4)n+2s)$ complete pairings. Assembling these factors and summing over vertex-set partitions via the exponential formula yields the stated expression. For $d=3$, we have $b=0$, which forces $n_0=0$, as a cubic vertex cannot accommodate four distinct points. (In exact arithmetic, Theorem~\ref{thm:secondmoment} agrees integer-for-integer with brute-force enumeration of all $M(12)=10{,}395$ pairings for $(n,d)=(4,3)$ and $(3,4)$ at $u\in\{\tfrac12,1,2\}$.)
\end{proof}

Write $M(2m)=(2m)!/(m!\,2^m)$ and $\E Z(u)=\sum_P Z(u,P)/M(dn)$. From Theorem~\ref{thm:moments},
\begin{equation}\label{eq:EZ-u}
\E Z(u)\;=\; n!\,[x^n]\,(1-x)^{-u/2}\,(d(d-1))^n\,\frac{M((d-2)n)}{M(dn)}
\end{equation}
in the multigraph (pairing) model; the girth-$\ge3$ version replaces $(1-x)^{-u/2}$ with $e^{-u(x/2+x^2/4)}(1-x)^{-u/2}$ as in Theorem~\ref{thm:girth}. Lemma~\ref{lem:delta} therefore predicts the small-subgraph-conditioning second-moment ratio
\begin{equation}\label{eq:Cu}
\frac{\E\,Z(u)^2}{\bigl(\E\,Z(u)\bigr)^2} \;\longrightarrow\;
C(u) \;=\; \exp\Bigl(\sum_{i\ge1}\lambda_i\,\delta_i(u)^2\Bigr)
\;=\; \Bigl(\frac{d-1}{d-2}\Bigr)^{\!\frac{1+(u-1)^2}{2}}\Bigl(\frac{d-1}{d}\Bigr)^{\!u-1}.
\end{equation}
Two consistency checks follow immediately: $C(u)$ is the constant for the \emph{pairing model}, where the product runs over all $i\ge1$. At $u=1$, we have $C(1)=\sqrt{(d-1)/(d-2)}$, and removing the factors for $i=1,2$ (the effect of conditioning on simplicity) yields
\[
\exp\Bigl(\sum_{i\ge3}\lambda_i(\delta_i^{2\mathrm F})^2\Bigr)
\;=\;
e^{-(2d-1)/(4(d-1)^2)}\sqrt{(d-1)/(d-2)},
\]
recovering the second-moment constant of Theorem~1 of Robalewska for simple $G(n,d)$ \cite{robalewska1996}, whose Theorem~3 accordingly features the product over $i\ge3$. As $u\to0^+$, $C(0^+)=d/(d-2)$ is the pairing Hamilton-cycle constant; the corresponding value for simple graphs is $C(0^+)e^{-\lambda_1(\delta_1^{\mathrm H})^2}=\tfrac{d}{d-2}e^{-2/(d-1)}$, cf.\ \cite{frieze1996,robinson_wormald1994}.

\begin{proposition}\label{prop:second-moment-asymp}
Fix $d\ge 3$. There exists $\varepsilon=\varepsilon(d)>0$ such that, in the pairing model,
\[
R_n(u)\;:=\;\frac{\E Z(u)^2}{\bigl(\E Z(u)\bigr)^2}\;\longrightarrow\; C(u)
\]
as $n\to\infty$, uniformly for $u$ in the closed disk $\overline{D}(1,\varepsilon)=\{u\in\mathbb C:|u-1|\le\varepsilon\}$, where $C(u)$ is the entire function \eqref{eq:Cu}. In particular, the convergence holds uniformly on the real interval $[1-\varepsilon,1+\varepsilon]$.
\end{proposition}

The full proof of Proposition~\ref{prop:second-moment-asymp} is given in Appendix~\ref{sec:laplace}; we summarise the main steps here. In terms of the rescaled variables $\alpha=n_0/n$ and $\beta=p/n$ (the proportions of free vertices and purple paths), Stirling's formula shows that each interior summand of $R_n(u)$ has the form $n^{-1}\exp\bigl(n\,\eta(\alpha,\beta)+\Psi_u(\alpha,\beta)+o(1)\bigr)$, where $\sigma=1-\alpha-\beta$ and $\hat\gamma=1-\alpha-2\beta$. The leading exponent $\eta$ does not depend on $u$, because $u$ enters only through shifts in the arguments of Gamma functions and hence only through the $O(1)$ term $\Psi_u$. As a result, $\eta$ attains its unique global maximum on the domain of summation at
\[
\alpha^*\;=\;\frac{(d-2)(d-3)}{d(d-1)},\qquad
\beta^*\;=\;\frac{2(d-2)}{d(d-1)},\qquad
\sigma^*\;=\;\frac{2}{d},
\]
with negative-definite Hessian $H=D^2\eta(\alpha^*,\beta^*)$ independent of $u$, and the tail bounds outside a neighbourhood of $(\alpha^*,\beta^*)$ are likewise uniform in $u$, while $u\mapsto\Psi_u(\alpha^*,\beta^*)$ is holomorphic. The maximiser $(\alpha^*,\beta^*)$ coincides with the point $c=(d-2)(d-3)n/(d(d-1))$, $k=2(d-2)n/(d(d-1))$ about which Robalewska expands her second-moment sum \cite{robalewska1996}, and the Gaussian integral at $(\alpha^*,\beta^*)$ evaluates to $C(u)$. (For $d=3$, where $b=0$ forces $n_0=0$, the sum is one-dimensional with $\beta^*=1/3$.)

\begin{lemma}\label{lem:mixed-second}
Fix $d\ge 3$ and let $\varepsilon=\varepsilon(d)>0$ be as in Proposition~\ref{prop:second-moment-asymp}. For $u,v\in\overline{D}(1,\varepsilon)$,
\begin{equation}\label{eq:mixed-R}
R_n(u,v)\;:=\;\frac{\E\bigl[Z(u)Z(v)\bigr]}{\E Z(u)\,\E Z(v)}
\;\longrightarrow\;
\exp\Bigl(\sum_{i\ge1}\lambda_i\,\delta_i(u)\,\delta_i(v)\Bigr)
\end{equation}
uniformly on the product of disks. In particular, for every $u\in\overline{D}(1,\varepsilon)$,
\begin{equation}\label{eq:mod-second}
\frac{\E\bigl|Z(u)\bigr|^2}{\bigl|\E Z(u)\bigr|^2}
\;\longrightarrow\;
\exp\Bigl(\sum_{i\ge1}\lambda_i\,\bigl|\delta_i(u)\bigr|^2\Bigr).
\end{equation}
\end{lemma}

\begin{proof}
The proof of Theorem~\ref{thm:secondmoment} adapts verbatim to the bilinear weight $u^{X(F_1)}v^{X(F_2)}$: purple cycles contribute $uv$ instead of $u^2$, and the two monochromatic completions contribute Pochhammer factors $(\tfrac u2)_{n_0+p}$ and $(\tfrac v2)_{n_0+p}$ instead of a square. The exponent $\eta$ remains unchanged, while the $O(1)$ term becomes
\begin{equation}\label{eq:Psi-uv}
\Psi_{u,v}(\alpha,\beta)
\;=\;
-\log(2\pi)-\tfrac12\log(\alpha\beta\hat\gamma)
+\tfrac{u+v-2}{2}\log(\alpha+\beta)
+\tfrac{uv-1}{2}\log(\sigma/\beta),
\end{equation}
where $\alpha,\beta,\sigma,\hat\gamma$ are the proportions defined in Appendix~\ref{sec:laplace}, cf.\ \eqref{eq:Psi-u}; $\Psi_{u,v}$ is jointly holomorphic for $|u-1|,|v-1|\le\varepsilon$ and coincides with $\Psi_u$ on the diagonal $v=u$. At the maximiser \eqref{eq:maximizer}, one has $\alpha^*+\beta^*=(d-2)/d$ and $\sigma^*/\beta^*=(d-1)/(d-2)$, so
\begin{align*}
e^{\Psi_{u,v}(\alpha^*,\beta^*)}\,\frac{2\pi}{\sqrt{\det(-H)}}
&\;=\;
C(1)\,\Bigl(\frac{d-2}{d}\Bigr)^{(u+v-2)/2}
\Bigl(\frac{d-1}{d-2}\Bigr)^{(uv-1)/2}\\
&\;=\;
\exp\Bigl(\sum_{i\ge1}\lambda_i\,\delta_i(u)\,\delta_i(v)\Bigr),
\end{align*}
which proves \eqref{eq:mixed-R}. (The second equality is the polarisation of the identity $e^{\Psi_u^*}\,2\pi/\sqrt{\det(-H)}=C(u)$ established in Appendix~\ref{sec:laplace}; for $d=3$, the same evaluation holds with the one-dimensional leading factor established there.)

For \eqref{eq:mod-second}, set $v=\overline{u}$. Since $\delta_i(\overline{u})=\overline{\delta_i(u)}$, the exponent becomes $\sum_i\lambda_i|\delta_i(u)|^2$. Furthermore, $Z(\overline{u})=\overline{Z(u)}$ (the cycle counts being real), which yields $\E Z(u)Z(\overline{u})=\E|Z(u)|^2$.
\end{proof}

\section{The conditional distribution of $X$ given $G$}\label{sec:typical}

We now pass from a uniform pair $(P,F)$ to the conditional distribution of $X$ given $G$, that is, to a uniformly random 2-factor of a fixed graph $G$, where $G$ is itself a uniformly random simple $d$-regular graph. The argument uses small subgraph conditioning \cite{janson1995,wormald1999}, following Robalewska's analysis of 2-factors \cite{robalewska1996}. We first determine the conditional probability of Hamiltonicity (Section~\ref{sec:typical-ham}) and establish small subgraph conditioning for the generating polynomial $Z(u)$ at real $u$ near $1$. In Section~\ref{sec:ssc-poly} we extend this convergence to a complex disk around $u=1$ by a normal-family compactness argument; differentiating at $u=1$ via Cauchy's estimate then gives the conditional mean, and evaluating at $u=e^{it/\sigma_n}$ gives the conditional central limit theorem (Section~\ref{sec:typical-clt}). The required second moment is the one computed in Section~\ref{sec:variance}.

\subsection{The conditional probability of a Hamilton cycle}\label{sec:typical-ham}

\begin{proposition}\label{prop:typical-ham}
Fix $d\ge3$, and let $G$ be a uniformly random simple $d$-regular graph on $n$ vertices. Then, as $n\to\infty$,
\[
\sqrt n\;\frac{Z_1(G)}{Z(G)} \;\Longrightarrow\; \Lambda
\;=\; \frac{\sqrt\pi\,e^{3/4}}{2}\prod_{i\ge3}\left(\frac{1+\delta_i^{\mathrm H}}{1+\delta_i^{2\mathrm F}}\right)^{\!\xi_i} e^{\lambda_i(\delta_i^{2\mathrm F}-\delta_i^{\mathrm H})},
\]
a strictly positive, almost surely finite random variable, where $(\xi_i)_{i\ge3}$ are independent $\mathrm{Poisson}(\lambda_i)$ random variables (the distributional limits of the short-cycle counts $(X_i)_{i\ge3}$) and $\delta_i^{\mathrm H}$, $\delta_i^{2\mathrm F}$ are given by Lemma~\ref{lem:delta}. In particular, for a uniformly random 2-factor $F$ of $G$,
\[
\Prob(F\text{ Hamiltonian}\mid G) \;=\; \Theta_p(n^{-1/2}).
\]
\end{proposition}

\begin{proof}
By small subgraph conditioning in the pairing model, $Z_1/\E Z_1\Rightarrow W_1$ (Robinson--Wormald \cite{robinson_wormald1994}) and $Z/\E Z\Rightarrow W$ (Robalewska \cite{robalewska1996}; re-proved below as the $u=1$ case of Lemma~\ref{lem:ssc-Zu}), both jointly with the short-cycle counts $(X_i)$ in the sense of Janson \cite[Theorem~1]{janson1995}, where $W_Y=\prod_{i\ge1}(1+\delta_i)^{X_i}e^{-\lambda_i\delta_i}$ is evaluated at the cycle counts of the same graph. Hence
\[
\frac{Z_1}{Z} = \frac{\E Z_1}{\E Z}\cdot\frac{W_1}{W}\,(1+o_p(1)),
\]
with expectations taken under the pairing model, and $\E Z_1/\E Z=\Prob(X=1)=\tfrac12\sqrt{\pi/n}\,(1+o(1))$ by Corollary~\ref{cor:ham}. Conditioning on simplicity corresponds to conditioning on $X_1=X_2=0$, an event of limiting probability $e^{-\lambda_1-\lambda_2}>0$; thus, the $o_p(1)$ term is preserved (Lemma~\ref{lem:simple-transfer}), and the factors for $i=1,2$ in $W_1/W$ become deterministic:
\[
\frac{e^{-\lambda_1\delta_1^{\mathrm H}-\lambda_2\delta_2^{\mathrm H}}}{e^{-\lambda_1\delta_1^{2\mathrm F}-\lambda_2\delta_2^{2\mathrm F}}}
\;=\;\frac{e^{1}}{e^{1/2}\,e^{-1/4}}\;=\;e^{3/4},
\]
using $\lambda_1\delta_1^{\mathrm H}=-1$, $\delta_2^{\mathrm H}=0$, $\lambda_1\delta_1^{2\mathrm F}=-\tfrac12$, and $\lambda_2\delta_2^{2\mathrm F}=\tfrac14$. This converts $\tfrac12\sqrt{\pi/n}$ into the girth-$\ge3$ constant $\tfrac{\sqrt\pi e^{3/4}}{2}n^{-1/2}$ of Corollary~\ref{cor:ham} and leaves the product over $i\ge3$. Since $(X_i)_{i\ge3}$ converges in distribution to $(\xi_i)_{i\ge3}$, that product converges in distribution to the product in the definition of $\Lambda$. To see that $\Lambda\in(0,\infty)$ almost surely, note that $1+\delta_i^{\mathrm H}>0$ and $1+\delta_i^{2\mathrm F}>0$ for all $i\ge3$ and $d\ge3$. With $\rho_i = \log\bigl((1+\delta_i^{\mathrm H})/(1+\delta_i^{2\mathrm F})\bigr) = -\,(d-1)^{-i} + O((d-1)^{-2i})$, the summands of $\sum_{i\ge3} \bigl(\xi_i\rho_i - \lambda_i(\delta_i^{\mathrm H}-\delta_i^{2\mathrm F})\bigr)$ are independent with expectations $O\bigl(i^{-1}(d-1)^{-i}\bigr)$ and variances $\lambda_i\rho_i^2 = O\bigl(i^{-1}(d-1)^{-i}\bigr)$, so the series converges almost surely by Kolmogorov's two-series theorem.
\end{proof}

Note the cancellation implicit in $\Lambda$: individually, $\sum_i\lambda_i(\delta_i^{2\mathrm F}-\delta_i^{\mathrm H}) = \sum_i \tfrac{1}{2i} = \infty$, but each divergent centering term is cancelled by the expectation of the corresponding $\xi_i$-term, by the same mechanism that separates the multigraph and girth-$\ge3$ constants in Theorem~\ref{thm:girth}. What Proposition~\ref{prop:typical-ham} does \emph{not} establish is the conditional mean or the central limit theorem for the full cycle count; to obtain those, we apply small subgraph conditioning to the generating polynomial $Z(u)=\sum_c u^cZ_c(G)$ for $u$ near $1$, whose first moment is given by Theorem~\ref{thm:pgf} and whose second moment is given by Theorem~\ref{thm:secondmoment}.

\subsection{Small subgraph conditioning for $Z(u)$ on a complex disk}\label{sec:ssc-poly}

\begin{lemma}\label{lem:ssc-Zu}
Fix $d\ge 3$ and let $u\in[1-\varepsilon,1+\varepsilon]$ be as in Proposition~\ref{prop:second-moment-asymp}. Let $(\xi_i)_{i\ge1}$ be independent $\mathrm{Poisson}(\lambda_i)$ variables (the distributional limits of the short-cycle counts $(X_i)_{i\ge1}$), and define
\[
W(u)\;=\;W(u;\xi)\;=\; \prod_{i\ge 1}\bigl(1+\delta_i(u)\bigr)^{\xi_i}\,e^{-\lambda_i\delta_i(u)} .
\]
Then in the pairing model, jointly in distribution,
\[
\Bigl(\frac{Z(u)}{\E Z(u)},\,(X_i)_{i\ge1}\Bigr)\;\Longrightarrow\;\bigl(W(u),\,(\xi_i)_{i\ge1}\bigr).
\]
In particular, $Z(u)/\E Z(u)\Rightarrow W(u)$ in distribution, with $\E W(u)=1$ and $W(u)>0$ almost surely. The family $\{W(u)\}$ is defined jointly on the sequence of independent Poisson variables $(\xi_i)_{i\ge1}$, and the map $u\mapsto W(u)$ is almost surely $C^1$ on the interval.
\end{lemma}

\begin{proof}
Proposition~\ref{prop:second-moment-asymp} gives
\[
\frac{\E Z(u)^2}{(\E Z(u))^2}\;\longrightarrow\; C(u)=\exp\Bigl(\sum_i\lambda_i\delta_i(u)^2\Bigr),
\]
Lemma~\ref{lem:delta} gives $\E[Z(u)X_i]/\E Z(u)\to\lambda_i(1+\delta_i(u))$, and Lemma~\ref{lem:factorial} provides the complete factorial-moment hierarchy
\[
\frac{\E\bigl[Z(u)\,(X_i)_r\bigr]}{\E Z(u)}\;\longrightarrow\;\lambda_i^r\bigl(1+\delta_i(u)\bigr)^r
\]
together with the joint moments across finitely many cycle lengths. Since $u>0$, the quantity $Y=Z(u)=\sum_c u^cZ_c$ is nonnegative, $\delta_i(u)>-1$ for all $i$, and $\sum_i\lambda_i\delta_i(u)^2<\infty$; these are precisely hypotheses (A1)--(A4) of Janson's small subgraph conditioning theorem \cite[Theorem~1]{janson1995}, whose conclusion is the stated joint convergence.

To verify the regularity of $W(u)$, note that $\delta_i(u)=(u-1+(-1)^i)/(d-1)^i$ is affine in $u$, and for $u$ near $1$ one has $1+\delta_i(u)>0$ uniformly in $i\ge 1$ for $d\ge 3$. We may therefore write $\log W(u)=\sum_i\bigl[(\xi_i-\lambda_i)\delta_i(u)+\xi_i\bigl(\log(1+\delta_i(u))-\delta_i(u)\bigr)\bigr]$. The first series decomposes as $(u-1)\sum_i(\xi_i-\lambda_i)(d-1)^{-i}+\sum_i(\xi_i-\lambda_i)(-1)^i(d-1)^{-i}$, which are two almost surely convergent random series (with variances $\sum_i\lambda_i(d-1)^{-2i}<\infty$) multiplied by entire functions of $u$; the second series converges absolutely and locally uniformly in $u$, since $|\log(1+\delta)-\delta|\le|\delta|^2$ and $\E\sum_i\xi_i|\delta_i(u)|^2<\infty$. Hence, $u\mapsto\log W(u)$ is almost surely real-analytic near $1$, and
\[
\partial_u\log W(u)
\;=\;
\sum_{i\ge 1}\Biggl(\xi_i\frac{\delta_i'(u)}{1+\delta_i(u)}-\lambda_i\delta_i'(u)\Biggr),
\qquad \delta_i'(u)=(d-1)^{-i},
\]
which is an almost surely convergent series, as for $\Lambda$ in Proposition~\ref{prop:typical-ham}.
\end{proof}

We next extend Lemma~\ref{lem:ssc-Zu} to uniform convergence on a complex disk centered at $u=1$ and deduce Theorems~\ref{thm:typical-mean} and~\ref{thm:typical-clt}. Small subgraph conditioning itself is invoked only for real $u>0$, where $Z(u)$ is nonnegative; the passage to complex $u$ is a compactness argument for random holomorphic functions, whose only additional input is the uniform bound on the complex second moment provided by Lemma~\ref{lem:mixed-second}. The conditional mean then follows from Cauchy's estimate for the derivative at $u=1$, and the CLT requires only $u_n=e^{it/\sigma_n}\to1$.

\begin{lemma}\label{lem:ssc-Zu-complex}
Fix $d\ge 3$, and choose $\varepsilon\in(0,1)$ small enough that the conclusion of Proposition~\ref{prop:second-moment-asymp} holds and $1+\delta_i(u)\ne0$ for all $i\ge1$ and $|u-1|\le\varepsilon$. Let $D'=D(1,\varepsilon/2)$ and let $\mathcal H(D')$ denote the space of holomorphic functions on $D'$ endowed with the topology of locally uniform convergence. Set $\widehat Y_n(u)=Z(u)/\E Z(u)$. With $(\xi_i)$ and $W(u)=W(u;\xi)$ as in Lemma~\ref{lem:ssc-Zu}, the function $u\mapsto W(u)$ extends almost surely to a zero-free holomorphic function on $D'$ satisfying $W(1)>0$, and in the pairing model we have
\begin{equation}\label{eq:complex-ssc}
\bigl(\widehat Y_n,\,(X_i)_{i\ge1}\bigr)\;\Longrightarrow\;\bigl(W,\,(\xi_i)_{i\ge1}\bigr)
\qquad\text{in distribution on }\mathcal H(D')\times\mathbb Z_{\ge0}^{\mathbb N}.
\end{equation}
In particular, $\sup_{|u-1|\le\varepsilon/4}|\widehat Y_n(u)|$ and $\sup_{|u-1|\le\varepsilon/4}|\widehat Y_n'(u)|$ are $O_p(1)$, $\widehat Y_n(1)\Rightarrow W(1)>0$, and $\widehat Y_n(u_n)-\widehat Y_n(1)=O_p(|u_n-1|)$ for every deterministic sequence $u_n\to1$.
\end{lemma}

\begin{proof}
Note first that since $|\delta_i(u)-\delta_i(1)|\le|u-1|(d-1)^{-i}$ and $1+\delta_i(1)\ge 1-(d-1)^{-1}>0$, choosing $\varepsilon>0$ sufficiently small ensures $1+\delta_i(u)\ne0$ for all $i\ge1$ and $|u-1|\le\varepsilon$.

\textit{The limit.} As in the proof of Lemma~\ref{lem:ssc-Zu}, $\log W(u)=(u-1)A+B+\sum_i\xi_i\bigl(\log(1+\delta_i(u))-\delta_i(u)\bigr)$, where $A=\sum_i(\xi_i-\lambda_i)(d-1)^{-i}$ and $B=\sum_i(\xi_i-\lambda_i)(-1)^i(d-1)^{-i}$ converge almost surely, and the remaining series converges absolutely and uniformly on $\overline D(1,\varepsilon)$ (using principal logarithms; $|\delta_i(u)|\le(1+\varepsilon)(d-1)^{-i}\le\tfrac12$ for $i\ge2$ after shrinking $\varepsilon$, and the term $i=1$ is a single holomorphic function). Hence, $W=e^{\log W}$ is almost surely holomorphic and zero-free on $D'$, and $W(1)>0$.

\textit{Tightness.} By Lemma~\ref{lem:mixed-second}, $\sup_{|u-1|\le\varepsilon}\E|\widehat Y_n(u)|^2\le M$ for all sufficiently large $n$, where $M=M(d,\varepsilon)<\infty$. For $z\in D'$, the disk $D(z,\varepsilon/2)$ is contained in $\overline D(1,\varepsilon)$, so the mean-value property of the holomorphic function $\widehat Y_n$ yields
\begin{align*}
|\widehat Y_n(z)|
&\;\le\;
\frac{4}{\pi\varepsilon^2}\int_{D(z,\varepsilon/2)}|\widehat Y_n|\,dA
\;\le\;
\frac{4}{\pi\varepsilon^2}\int_{\overline D(1,\varepsilon)}|\widehat Y_n|\,dA,\\
\text{hence}\qquad
\E\sup_{z\in D'}|\widehat Y_n(z)|
&\;\le\;
4\sqrt M .
\end{align*}
By Markov's inequality and Montel's theorem (a locally bounded family is relatively compact in $\mathcal H(D')$), the sequence of distributions of $\widehat Y_n$ is tight in $\mathcal H(D')$; the distributions of $(X_i)_{i\ge1}$ are tight in $\mathbb Z_{\ge0}^{\mathbb N}$ since each coordinate converges. Thus, by Prohorov's theorem, every subsequence of $(\widehat Y_n,(X_i))$ possesses a further subsequence converging in distribution to some limit $(\widehat Y_\infty,(\xi_i))$, where the variables $(\xi_i)$ are independent $\mathrm{Poisson}(\lambda_i)$ by the classical short-cycle limit.

\textit{Identification of the limit.} Fix a real $u\in(1-\varepsilon/2,1+\varepsilon/2)$. Evaluation at $u$ is a continuous functional on $\mathcal H(D')$, so along the chosen subsequence, $(\widehat Y_n(u),(X_i))\Rightarrow(\widehat Y_\infty(u),(\xi_i))$. On the other hand, Lemma~\ref{lem:ssc-Zu} asserts that $(\widehat Y_n(u),(X_i))\Rightarrow(W(u;\xi),(\xi_i))$ along the full sequence. The two limit distributions therefore coincide, and $W(u;\xi)$ is a measurable function of $(\xi_i)$; applying the equality of distributions to the bounded functional $\min\{1,|y-W(u;x)|\}$ of $(y,x)$ gives $\E\min\{1,|\widehat Y_\infty(u)-W(u;\xi)|\}=0$, that is, $\widehat Y_\infty(u)=W(u;\xi)$ almost surely. Taking $u$ in a countable dense subset of the real interval, the almost surely holomorphic function $\widehat Y_\infty-W(\cdot\,;\xi)$ vanishes on a set having an accumulation point in $D'$, and thus vanishes identically on $D'$ by the identity theorem. Consequently, every subsequential limit is $(W,(\xi_i))$, which establishes \eqref{eq:complex-ssc}.

\textit{Consequences.} The functionals $f\mapsto\sup_{|u-1|\le\varepsilon/4}|f|$, $f\mapsto\sup_{|u-1|\le\varepsilon/4}|f'|$, and $f\mapsto f(1)$ are continuous on $\mathcal H(D')$ (by Cauchy's estimate for the second); hence, by the continuous mapping theorem, these three quantities converge in distribution to the corresponding almost surely finite functionals of $W$. Convergence in distribution to a finite limit implies tightness, that is, boundedness in probability $O_p(1)$. Finally, $|\widehat Y_n(u_n)-\widehat Y_n(1)|\le|u_n-1|\sup_{|u-1|\le\varepsilon/4}|\widehat Y_n'(u)|$ as soon as $|u_n-1|\le\varepsilon/4$.
\end{proof}

\subsection{Conditional mean, central limit theorem, and variance}\label{sec:typical-clt}

\begin{theorem}\label{thm:typical-mean}
Fix $d\ge 3$. Let $G$ be a uniformly random simple $d$-regular graph on $n$ vertices, and let $X$ be the number of cycles of a uniformly random 2-factor of $G$. With $\mu_n^{\mathrm{pair}}=\tfrac12\log n+\tfrac{\gamma}{2}+\log 2$ as in Corollary~\ref{cor:quasipower-input}, write
\[
\mu_n\;=\;\mu_n^{\mathrm{pair}}-\tfrac34\;=\;\tfrac12\log n+\tfrac{\gamma}{2}+\log 2-\tfrac34.
\]
Let $(\xi_i)_{i\ge3}$ be independent $\mathrm{Poisson}(\lambda_i)$ variables, the distributional limits of $(X_i)_{i\ge3}$, and let $W_{\ge3}(u)=\prod_{i\ge3}(1+\delta_i(u))^{\xi_i}e^{-\lambda_i\delta_i(u)}$ be the product of Lemma~\ref{lem:ssc-Zu} truncated at $i\ge3$. Then
\[
\E[X\mid G]-\mu_n
\;\Longrightarrow\;
\partial_u\log W_{\ge3}(u)\Big|_{u=1}
\qquad\text{in distribution.}
\]
The right-hand side equals
\[
\sum_{i\ge3}\Bigl(\xi_i\,\frac{(d-1)^{-i}}{1+(-1)^i(d-1)^{-i}}-\frac1{2i}\Bigr).
\]
The series converges almost surely and has finite variance
\begin{equation}\label{eq:varW}
\begin{aligned}
\sum_{i\ge 3}\lambda_i\,\frac{\delta_i'(1)^2}{(1+\delta_i(1))^2}
&\;=\;
\sum_{i\ge 3}\frac{(d-1)^{-i}}{2i}\,\frac{1}{\bigl(1+(-1)^i(d-1)^{-i}\bigr)^2}\\
&\;<\;\infty.
\end{aligned}
\end{equation}
In particular, $\E[X\mid G]=\tfrac12\log n+\tfrac\gamma2+\log2-\tfrac34+O_p(1)$: for every $\epsilon>0$, there exists $C=C(d,\epsilon)<\infty$ such that
\begin{equation}\label{eq:typical-mean-ub}
\Prob\bigl(\,|\E[X\mid G]-\mu_n|\le C\bigr)\;\ge\;1-\epsilon
\end{equation}
for all sufficiently large $n$.
\end{theorem}

\begin{proof}
Under the conditional distribution of a uniform 2-factor of the fixed graph $G$,
\begin{align*}
\E\bigl[u^{X}\bigm|G\bigr]
&\;=\;\frac{Z(u)}{Z(1)},\\
\E[X\mid G]
&\;=\;\partial_u\log Z(u)\Big|_{u=1}
\;=\;\partial_u\log\E Z(u)\Big|_{u=1}+\frac{\widehat Y_n'(1)}{\widehat Y_n(1)},
\end{align*}
with $\widehat Y_n(u)=Z(u)/\E Z(u)$ as defined in Lemma~\ref{lem:ssc-Zu-complex}; the final identity follows from $\log Z(u)=\log\E Z(u)+\log\widehat Y_n(u)$, valid near $u=1$ where $\widehat Y_n(1)=Z/\E Z>0$. In the pairing model, $\E Z(u)/\E Z(1)=\E u^X$ is the generating function given by Theorem~\ref{thm:pgf}, so that
\begin{equation}\label{eq:pair-deriv}
\partial_u\log\E Z(u)\Big|_{u=1}\;=\;\E X\;=\;\mu_n^{\mathrm{pair}}+o(1).
\end{equation}

\textit{Pairing model.} The map $\Phi(f)=f'(1)/f(1)$ is continuous on the open subset $\{f(1)\ne0\}$ of $\mathcal H(D')$, and $W(1)>0$ almost surely. By \eqref{eq:complex-ssc} and the continuous mapping theorem,
\begin{align*}
&\Bigl(\E[X\mid G]-\mu_n^{\mathrm{pair}},\,(X_i)_{i\ge1}\Bigr)\\
&\qquad\Longrightarrow\;
\Bigl(\partial_u\log W(u)\big|_{u=1},\,(\xi_i)_{i\ge1}\Bigr),\\
&\partial_u\log W(1)
\;=\;
\sum_{i\ge1}\Bigl(\xi_i\frac{\delta_i'(1)}{1+\delta_i(1)}-\lambda_i\delta_i'(1)\Bigr).
\end{align*}

\textit{Simple graphs.} Conditioning the pairing on simplicity corresponds to conditioning on the event $\{X_1=X_2=0\}$, which is determined by two discrete coordinates of the second component and has limiting probability $e^{-\lambda_1-\lambda_2}>0$. The joint convergence in distribution therefore transfers to the conditional distributions: conditionally on simplicity, $\E[X\mid G]-\mu_n^{\mathrm{pair}}$ converges in distribution to $\partial_u\log W(1)$ conditioned on $\xi_1=\xi_2=0$. On this event, the terms for $i=1,2$ are deterministic,
\[
-\lambda_1\delta_1'(1)-\lambda_2\delta_2'(1)\;=\;-\tfrac12-\tfrac14\;=\;-\tfrac34,
\]
and the remaining sum is $\partial_u\log W_{\ge3}(1)$, whose summands involve only $(\xi_i)_{i\ge3}$ and are independent of $(\xi_1,\xi_2)$. Since $\mu_n=\mu_n^{\mathrm{pair}}-\tfrac34$, this establishes the convergence in distribution.

\textit{The limit.} The series converges almost surely and has variance \eqref{eq:varW} by the calculation in Lemma~\ref{lem:ssc-Zu}, as the summands are independent with $\operatorname{Var}(\xi_i)=\lambda_i$. Convergence in distribution to an almost surely finite limit implies tightness, which yields \eqref{eq:typical-mean-ub} with $C=C(d,\epsilon)$ chosen as any constant exceeding the $(1-\epsilon/2)$-quantile of the absolute value of the series. Because $\xi_3$ is unbounded, the series is likewise unbounded, so no deterministic constant $C$ can hold with probability $1-o(1)$: the fluctuation of $\E[X\mid G]$ about $\mu_n$ is of order one and is governed by the short cycles of $G$.
\end{proof}

\begin{theorem}\label{thm:typical-clt}
Fix $d\ge3$ and let $G$ be a uniformly random \emph{simple} $d$-regular graph on $n$ vertices. If $X$ is the number of cycles of a uniformly random 2-factor of $G$, then
\[
\frac{X-\mu_n}{\sigma_n}\;\Bigm|\;G
\;\Longrightarrow\; N(0,1)
\]
in probability (that is, for every bounded continuous function $h$, $\E[h((X-\mu_n)/\sigma_n)\mid G]\to\E h(N)$ in probability over $G$, where $N\sim N(0,1)$), with
\[
\mu_n=\tfrac12\log n+\tfrac{\gamma}{2}+\log 2-\tfrac34,\qquad
\sigma_n=\sqrt{\tfrac12\log n}.
\]
\end{theorem}

\begin{proof}
We work in the pairing model with the exact centering $m_n=\E X$ for a uniform pair from Corollary~\ref{cor:quasipower-input}; the transfer to simple graphs $G(n,d)$ and to the centering $\mu_n$ is carried out at the end. Fix $t\in\mathbb R$ and set $u_n=\exp(it/\sigma_n)$, so that $|u_n-1|=O(1/\sigma_n)\to0$. The conditional characteristic function given $G$ is
\begin{equation}\label{eq:phi-G}
\begin{aligned}
\phi_G(t)
&\;:=\;
\E\Bigl[\exp\Bigl(\frac{it(X-m_n)}{\sigma_n}\Bigr)\Bigm|G\Bigr]\\
&\;=\;
e^{-itm_n/\sigma_n}\,\frac{Z(u_n)}{Z(1)}\\
&\;=\;
\Bigl[e^{-itm_n/\sigma_n}\,\frac{\E Z(u_n)}{\E Z(1)}\Bigr]
\cdot\frac{\widehat Y_n(u_n)}{\widehat Y_n(1)} .
\end{aligned}
\end{equation}
The bracketed factor is the characteristic function of $(X-m_n)/\sigma_n$ for a uniform pair; by \eqref{eq:EZ-complex} with $s=it/\sigma_n$, it equals
\[
\exp\Bigl(-\tfrac12 v_n t^2/\sigma_n^2+O\bigl(|t|^3\sigma_n^{-3}\log n\bigr)\Bigr)\;\longrightarrow\;e^{-t^2/2},
\]
since $v_n/\sigma_n^2\to1$ and $\sigma_n^{-3}\log n=O((\log n)^{-1/2})$. For the random factor, Lemma~\ref{lem:ssc-Zu-complex} yields $\widehat Y_n(u_n)-\widehat Y_n(1)=O_p(\sigma_n^{-1})$ and $\widehat Y_n(1)\Rightarrow W(1)>0$, whence $1/\widehat Y_n(1)=O_p(1)$ and
\[
\frac{\widehat Y_n(u_n)}{\widehat Y_n(1)}\;=\;1+O_p(\sigma_n^{-1})\;\longrightarrow\;1
\qquad\text{in probability.}
\]
Fix $T<\infty$. Uniformly for $|t|\le T$, \eqref{eq:EZ-complex} with $s=it/\sigma_n$ gives
\[
\Bigl[e^{-itm_n/\sigma_n}\,\frac{\E Z(e^{it/\sigma_n})}{\E Z(1)}\Bigr]
=
\exp\Bigl(-\tfrac12 v_n t^2/\sigma_n^2+O_T(\sigma_n^{-3}\log n)\Bigr),
\]
which converges uniformly to $e^{-t^2/2}$, since $v_n/\sigma_n^2\to1$ and $\sigma_n^{-3}\log n\to0$. For all large $n$ the points $e^{it/\sigma_n}$ with $|t|\le T$ lie in $|u-1|\le\varepsilon/4$, so Lemma~\ref{lem:ssc-Zu-complex} yields
\[
\sup_{|t|\le T}
\bigl|\widehat Y_n(e^{it/\sigma_n})-\widehat Y_n(1)\bigr|
\;\le\;
\sup_{|u-1|\le\varepsilon/4}|\widehat Y_n'(u)|\,
\sup_{|t|\le T}|e^{it/\sigma_n}-1|
\;=\;O_p(\sigma_n^{-1}).
\]
Together with $\widehat Y_n(1)\Rightarrow W(1)>0$, this gives
\[
\sup_{|t|\le T}
\left|\frac{\widehat Y_n(e^{it/\sigma_n})}{\widehat Y_n(1)}-1\right|
\;=\;o_p(1),
\]
and therefore
\[
\sup_{|t|\le T}\bigl|\phi_G(t)-e^{-t^2/2}\bigr|
\;\longrightarrow\;0
\qquad\text{in probability.}
\]
Along any subsequence, a diagonal extraction produces a further subsequence on which this convergence holds almost surely for every integer $T$. For almost every such graph sequence the conditional characteristic functions converge pointwise on $\mathbb R$ to $e^{-t^2/2}$, which is continuous at the origin, so L\'evy's continuity theorem yields
\[
\mathcal L\bigl((X-m_n)/\sigma_n\bigm|G\bigr)\Rightarrow N(0,1)
\]
almost surely along that further subsequence. Tightness is part of this conclusion. Since every subsequence has such a further subsequence, the conditional distributions converge weakly to $N(0,1)$ in probability along the full sequence.

\textit{Transfer.} Conditioning on simplicity preserves every $o_p(1)$ and $O_p(1)$ asymptotic statement above (Lemma~\ref{lem:simple-transfer}), and no short-cycle factors appear in \eqref{eq:phi-G} other than through $\widehat Y_n$, so no adjustments for $i=1,2$ are needed. Finally, $m_n-\mu_n=\tfrac34+o(1)=O(1)$, so replacing the centering sequence $m_n$ by $\mu_n$ alters $(X-m_n)/\sigma_n$ by an error of $O(\sigma_n^{-1})\to0$. (Similarly, centering at the conditional mean $\E[X\mid G]=\mu_n+O_p(1)$ of Theorem~\ref{thm:typical-mean} shifts the normalised variable by $O_p(\sigma_n^{-1})$ and does not affect the limiting distribution.)
\end{proof}

\begin{proposition}\label{rem:typical-var}
Under the hypotheses of Theorem~\ref{thm:typical-clt},
\[
\operatorname{Var}(X\mid G)\;=\;\sigma_n^2\,\bigl(1+o_p(1)\bigr).
\]
\end{proposition}

\begin{proof}
Work first in the pairing model. Let $P_n$ be the uniform distribution on pairings, and let $Q_n$ be the distribution on pairings with density $dQ_n/dP_n=\widehat Y_n(1)=Z(1)/\E Z(1)$. The proof of Theorem~\ref{thm:typical-clt} gives $(X-m_n)/\sigma_n\mid G\Rightarrow N(0,1)$ in $P_n$-probability. Set $V_n=\operatorname{Var}(X\mid G)/\sigma_n^2$.

\emph{Lower bound.} Fix $\varepsilon>0$ and set $f_K(x)=\min(x^2,K^2)$. Choose $K$ so that $\E f_K(N)\ge1-\varepsilon/2$, where $N\sim N(0,1)$. Then
\[
M_{n,K}
\;:=\;
\E\Bigl[f_K\Bigl(\frac{X-m_n}{\sigma_n}\Bigr)\Bigm|G\Bigr]
\;\longrightarrow\;
\E f_K(N)
\]
in $P_n$-probability. By the pairing-model step in the proof of Theorem~\ref{thm:typical-mean}, $\E[X\mid G]-m_n=O_p(1)$. Writing
\[
r_n
\;=\;
\Bigl(\frac{\E[X\mid G]-m_n}{\sigma_n}\Bigr)^2,
\]
one has $r_n=o_p(1)$ under $P_n$ and $\E\bigl[((X-m_n)/\sigma_n)^2\bigm|G\bigr]=V_n+r_n$. Since $f_K(x)\le x^2$, the event $\{V_n<1-\varepsilon\}$ is contained in $\{M_{n,K}<1-3\varepsilon/4\}\cup\{r_n>\varepsilon/4\}$, and both of those probabilities tend to zero under $P_n$. Thus $\Prob_{P_n}(V_n<1-\varepsilon)\to0$.

\emph{Contiguity.} By Lemma~\ref{lem:ssc-Zu}, $\widehat Y_n(1)\Rightarrow W(1)$ with $W(1)>0$ almost surely and $\E W(1)=1$. Since also $\E_{P_n}\widehat Y_n(1)=1$, these densities are uniformly integrable, so every event of vanishing $P_n$-probability has vanishing $Q_n$-probability. Conversely, if $Q_n(A_n)\to0$, then for every $\eta>0$
\[
P_n(A_n)
\;\le\;
P_n\bigl(\widehat Y_n(1)\le\eta\bigr)+\eta^{-1}Q_n(A_n);
\]
letting $n\to\infty$ and then $\eta\downarrow0$ gives $P_n(A_n)\to0$, because $\widehat Y_n(1)\Rightarrow W(1)>0$. The two distributions are mutually contiguous.

\emph{Upper bound.} The uniform pair of Corollary~\ref{cor:clt} has graph marginal $Q_n$, and the conditional distribution of $X$ given $G$ is the one already used. The law of total variance therefore yields $\E_{Q_n}V_n\le\operatorname{Var}X/\sigma_n^2=1+O(1/\log n)$, where $\operatorname{Var}X$ is the variance for a uniform pair. Contiguity transfers the lower tail, so $\Prob_{Q_n}(V_n<1-\varepsilon)\to0$. Since $0\le(1-V_n)_+\le1$,
\[
\E_{Q_n}(V_n-1)_+
\;=\;
\E_{Q_n}(V_n-1)+\E_{Q_n}(1-V_n)_+
\;\le\;
O(1/\log n)+\varepsilon+o(1).
\]
As $\varepsilon>0$ is arbitrary, $\E_{Q_n}(V_n-1)_+\to0$, and Markov's inequality gives $\Prob_{Q_n}(V_n>1+\delta)\to0$ for every $\delta>0$. Reverse contiguity returns this upper tail to $P_n$. Combining the upper and lower tails under $P_n$ yields $V_n\to1$ in $P_n$-probability, and Lemma~\ref{lem:simple-transfer} transfers the conclusion to simple graphs.
\end{proof}

\section*{Acknowledgements}
We thank Guy Blachar and other colleagues for comments on earlier drafts of this paper.

\noindent\textbf{AI Usage Statement:} During the preparation of this manuscript, the authors used Google Gemini (Argon) to assist with language editing and clarity, xAI (Grok~4.6) to generate Figure~1 and the C++ code, and Anthropic Claude (Fable~5.2) to verify the mathematical proofs. The authors thoroughly reviewed, tested, and validated all AI-generated code, figures, and proof checks prior to incorporation and take full responsibility for the accuracy and contents of the paper.

\appendix

\section{Asymptotics of the second moment}\label{sec:laplace}

This appendix gives the full proof of Proposition~\ref{prop:second-moment-asymp}.

\begin{proof}[Proof of Proposition~\ref{prop:second-moment-asymp}]
We work throughout with the multigraph (pairing) expectations of \eqref{eq:EZ-u} and Theorem~\ref{thm:secondmoment}. Both $\E Z(u)$ and $\sum_P Z(u,P)^2$ are entire in $u$ for each fixed $n$ (finite sums of monomials $u^{X(F)}$), so $R_n(u)$ is meromorphic; on a small disk about $u=1$ one has $\E Z(u)\ne0$ for large $n$ by the first-moment formula \eqref{eq:EZ-u}, hence $R_n$ is holomorphic there.

\textit{Normalisation.} Write $N(dn)=M(dn)$. Then
\[
R_n(u)
\;=\;
\frac{N(dn)\,\sum_P Z(u,P)^2}{\bigl(\sum_P Z(u,P)\bigr)^2}.
\]
Substituting Theorem~\ref{thm:secondmoment} and $\sum_P Z(u,P)=n!\,[x^n](1-x)^{-u/2}\,a^n\,M((d-2)n)$ with $a=d(d-1)$ yields an explicit sum over lattice points $(n_0,p)\in\mathbb Z_{\ge0}^2$. The summand vanishes unless $n-n_0-2p\ge0$, so the support is the triangle
\[
\Delta\;:=\;\bigl\{\alpha,\beta\ge0:\ \alpha+2\beta\le1\bigr\},\qquad \alpha:=n_0/n,\quad \beta:=p/n.
\]
Write $\sigma:=1-\alpha-\beta$, $\hat\gamma:=1-\alpha-2\beta$ (where the hat distinguishes this proportion from Euler's constant $\gamma$), $\tau:=d-2-2\alpha-2\beta$, and $N:=n_0+p$, $m:=n-n_0$, $g:=m-2p=\hat\gamma n$, $s:=\sigma n$. The purple coefficient is exact: $[y^g](1-y)^{-\nu}=(\nu)_g/g!$ is a polynomial identity in $\nu$, hence
\[
[x^{m}]\,e^{(u^2/2)L(ax)}\,\hat\pi(x)^{p}
\;=\;
a^{m}\Bigl(\frac{c_e^2}{2a^2}\Bigr)^{p}\frac{(u^2/2+p)_{g}}{g!},
\]
and $4^p(c_e^2/(2a^2))^p=(2(d-2)^2)^p$. Combined with $(\tfrac u2)_N=\Gamma(N+u/2)/\Gamma(u/2)$ and $[x^n](1-x)^{-u/2}=(\tfrac u2)_n/n!$, every summand of $R_n(u)$ is a ratio of Gamma values, factorials, and matching numbers $M(\tau n)\,M(dn)/M((d-2)n)^2$.

\textit{The exponent.} Stirling's formula $\log\Gamma(z)=(z-\tfrac12)\log z-z+\tfrac12\log(2\pi)+O(|z|^{-1})$, uniformly for $|\arg z|\le\pi-\delta$ and $|z|\to\infty$, applied to every Gamma, factorial, and matching factor, produces two pieces. The $n\log n$ exponents cancel identically (degree $d+1$ in the numerator and the denominator). The coefficient of $n$, $\eta$, is independent of $u$: every $u$ enters through a shift $\Gamma(N+u/2)$, $\Gamma(n+u/2)$, $\Gamma(s+u^2/2)$, or $\Gamma(p+u^2/2)$, which contributes only an $O(1)$ term after the $\log n$ pieces cancel. Explicitly, for $d\ge4$ and $(\alpha,\beta)$ in any fixed compact subset of the interior of $\Delta$,
\begin{equation}\label{eq:eta}
\begin{aligned}
\eta(\alpha,\beta)
&\;=\;
2(\alpha+\beta)\log(\alpha+\beta)
+\sigma\log\sigma
+\tfrac\tau2\log\tau
+\tfrac d2\log d\\
&\qquad
+\alpha\log b
+\beta\log\bigl(2(d-2)^2\bigr)\\
&\qquad
-\alpha\log\alpha
-2\beta\log\beta
-\hat\gamma\log\hat\gamma
-\log a
-(d-2)\log(d-2),
\end{aligned}
\end{equation}
with $b=(d-2)(d-3)$ and the convention $0\log0=0$ on the boundary, and
\begin{equation}\label{eq:Psi-u}
\Psi_u(\alpha,\beta)
\;=\;
-\log(2\pi)-\tfrac12\log(\alpha\beta\hat\gamma)
+(u-1)\log(\alpha+\beta)
+\tfrac{u^2-1}{2}\log(\sigma/\beta).
\end{equation}
The prefactors $(2\pi n)^{1/2}$ from Stirling's formula contribute $n^{1/2}$ in the numerator (from $n!$) and $n^{3/2}$ in the denominator (from $n_0!$, $p!$, and $g!$), leaving a factor of $n^{-1}$, the area of a single lattice cell in the $(\alpha,\beta)$-plane. Thus, for $(\alpha,\beta)$ in a fixed compact subset of the open triangle, the summand is
\begin{equation}\label{eq:rho-summand}
\rho_n(\alpha,\beta;u)
\;=\;
n^{-1}\exp\bigl(n\,\eta(\alpha,\beta)+\Psi_u(\alpha,\beta)+o(1)\bigr),
\end{equation}
uniformly for $(\alpha,\beta)$ in any fixed compact subset of the open triangle and $|u-1|\le\varepsilon$ (principal branches, $\mathrm{Re}\,u/2>0$). In particular $\eta$ does not depend on $u$, so the maximiser computed next is the same for every such $u$, including in a complex disk.

\textit{Consecutive ratios.} Write $T(n_0,p)$ for the $u=1$ summand of $\sum_P Z(1,P)^2$. The Pochhammer recurrence $(\tfrac12)_{k+1}/(\tfrac12)_k=k+\tfrac12$, the matching identity $M(2k-2)/M(2k)=1/(2k-1)$, and the exact coefficient ratio give
\begin{align}
\frac{T(n_0+1,p)}{T(n_0,p)}
&\;\longrightarrow\;
f_\alpha
\;:=\;
\frac{(\alpha+\beta)^2\,b\,\hat\gamma}{\alpha\,\sigma\,\tau},
\label{eq:ratio-alpha}\\
\frac{T(n_0,p+1)}{T(n_0,p)}
&\;\longrightarrow\;
f_\beta
\;:=\;
\frac{(\alpha+\beta)^2\,2(d-2)^2\,\hat\gamma^2}{\beta^2\,\sigma\,\tau},
\label{eq:ratio-beta}
\end{align}
as $n\to\infty$ with $(\alpha,\beta)$ in a compact subset of the open triangle. Differentiating \eqref{eq:eta} recovers $\nabla\eta=(\log f_\alpha,\log f_\beta)$. Setting both limits to $1$ yields
\begin{align}
(\alpha+\beta)^2\,b\,\hat\gamma
&\;=\;
\alpha\,\sigma\,\tau,
\label{eq:crit-alpha}\\
(\alpha+\beta)^2\,2(d-2)^2\,\hat\gamma^2
&\;=\;
\beta^2\,\sigma\,\tau.
\label{eq:crit-beta}
\end{align}
These are solved by
\begin{equation}\label{eq:maximizer}
\alpha^*\;=\;\frac{(d-2)(d-3)}{d(d-1)},\qquad
\beta^*\;=\;\frac{2(d-2)}{d(d-1)},\qquad
\sigma^*\;=\;\frac{2}{d}.
\end{equation}
At this point $\hat\gamma^*=2/(d(d-1))$ and $\tau^*=(d-2)^2/d$; direct substitution confirms that both sides of \eqref{eq:crit-alpha} equal $2(d-3)(d-2)^3/(d^3(d-1))$ and both sides of \eqref{eq:crit-beta} equal $8(d-2)^4/(d^4(d-1)^2)$. Dividing the two critical equations gives the relation $2(d-2)\alpha\hat\gamma=(d-3)\beta^2$. With $r=\beta/\alpha$ this reads $\alpha=2(d-2)/\bigl((d-3)r^2+2(d-2)(1+2r)\bigr)$, and substituting into \eqref{eq:crit-alpha} and clearing denominators leaves the cubic
\[
\bigl((d-3)r-2\bigr)\Bigl((d-3)(2d-5)\,r^2+2(d-3)(2d-5)\,r+2(d-2)(d-4)\Bigr)\;=\;0 .
\]
For $d\ge4$ the quadratic factor has nonnegative coefficients and positive leading coefficient, hence no root $r>0$; so $r=2/(d-3)$ is the only root with $\alpha,\beta>0$ and $\hat\gamma>0$, recovering \eqref{eq:maximizer} uniquely in the open triangle. On the open edges, $f_\alpha\to+\infty$ as $\alpha\to0^+$ with $\beta$ fixed in $(0,1/2)$, $f_\beta\to+\infty$ as $\beta\to0^+$, and $f_\beta\to0$ as $\hat\gamma\to0^+$. At the three corners, $\eta(0,0)=\tfrac d2\log d-\tfrac{d-2}2\log(d-2)-\log(d(d-1))<0$, and likewise $\eta(1,0)<0$, $\eta(0,\tfrac12)<0$. Thus no boundary point is a local maximum, and \eqref{eq:maximizer} is the unique global maximiser on the compact set $\Delta$. Substituting the coordinates into \eqref{eq:eta} and collecting coefficients of $\log(d-3)$, $\log(d-1)$, $\log d$, $\log(d-2)$ and $\log 2$ gives respectively $\alpha^*-\alpha^*$, $\alpha^*+2\beta^*+\hat\gamma^*-1$, $\tfrac d2-\tfrac{2(d-2)}{d}-\tfrac2d-\tfrac{(d-2)^2}{2d}$, $\tfrac{2(d-2)}{d}+\tfrac{(d-2)^2}{d}-(d-2)$ and $\tfrac2d-\beta^*-\hat\gamma^*$, each identically zero, so $\eta(\alpha^*,\beta^*)=0$.

\textit{Hessian.} Differentiating $(\log f_\alpha,\log f_\beta)$ at \eqref{eq:maximizer} produces $H:=D^2\eta(\alpha^*,\beta^*)$ with entries
\begin{align}
H_{11}
&\;=\;
-\frac{d\bigl(d^4-9d^3+28d^2-34d+16\bigr)}{2(d-3)(d-2)^2},
\nonumber\\
H_{12}
&\;=\;
\frac{d(-2d^3+11d^2-16d+8)}{2(d-2)^2},
\nonumber\\
H_{22}
&\;=\;
\frac{d(-4d^3+19d^2-26d+12)}{2(d-2)^2}
\label{eq:hessian-entries}
\end{align}
for $d\ge4$. In particular
\begin{equation}\label{eq:hessian-det}
\begin{aligned}
\det H
&\;=\;
\frac{d^3(d-1)^2}{4(d-3)(d-2)}
\;>\;0,\\
\operatorname{tr} H
&\;=\;
-\frac{d\bigl(5d^4-40d^3+111d^2-124d+52\bigr)}{2(d-3)(d-2)^2}
\;<\;0
\end{aligned}
\end{equation}
for every $d\ge4$: both quartics $d^4-9d^3+28d^2-34d+16$ and $5d^4-40d^3+111d^2-124d+52$ are positive on $[4,\infty)$ (each has all real critical points in $(0,4)$, and the endpoint values at $d=4$ are positive). Thus $H\prec0$.

\textit{The neighbourhood of the maximiser.} Fix a width $\omega_n=n^{-2/5}$ (so $\sqrt{(\log n)/n}\ll\omega_n\ll1$). Let $\mathcal L_n$ be the lattice points in $\Delta$ with $\|(\alpha,\beta)-(\alpha^*,\beta^*)\|\le\omega_n$, and write $R_n^{\mathrm{loc}}(u)=\sum_{\mathcal L_n}\rho_n(\alpha,\beta;u)$. For large $n$ this neighbourhood lies in a fixed compact subset of the open triangle, at a definite distance from $\partial\Delta$. Taylor-expand $\eta$ about the maximiser: with $x=\sqrt n\,(\alpha-\alpha^*,\beta-\beta^*)\in\mathbb R^2$,
\begin{equation}\label{eq:taylor-eta}
n\,\eta(\alpha,\beta)
\;=\;
\tfrac12\,x^{\mathsf T}H x
+O\bigl(|x|^3 n^{-1/2}\bigr)
\;=\;
\tfrac12\,x^{\mathsf T}H x
+O\bigl(n^{-1/5}\bigr),
\end{equation}
since $|x|\le n^{1/10}$ on this neighbourhood and therefore $|x|^3 n^{-1/2}\le n^{-1/5}$, while $\Psi_u(\alpha,\beta)=\Psi_u(\alpha^*,\beta^*)+O(\omega_n)$ uniformly for $|u-1|\le\varepsilon$. The lattice has spacing $1/n$ in $(\alpha,\beta)$, hence spacing $n^{-1/2}$ in the $x$-plane. Substituting \eqref{eq:rho-summand}--\eqref{eq:taylor-eta},
\begin{align}
R_n^{\mathrm{loc}}(u)
&\;=\;
e^{\Psi_u(\alpha^*,\beta^*)}\sum_{x}
\exp\Bigl(\tfrac12 x^{\mathsf T}H x+O(n^{-1/5})\Bigr)\,n^{-1}
\nonumber\\
&\;=\;
e^{\Psi_u(\alpha^*,\beta^*)}
\int_{\|x\|\le n^{1/10}}
\exp\bigl(\tfrac12 x^{\mathsf T}H x\bigr)\,dx
\;\bigl(1+O(n^{-1/5})\bigr),
\label{eq:local-riemann}
\end{align}
uniformly for $|u-1|\le\varepsilon$. Extending the integral to $\mathbb R^2$ costs $O(\exp(-c\,n^{1/5}))$ because $H\prec0$, and
\[
\int_{\mathbb R^2}\exp\bigl(\tfrac12 x^{\mathsf T}H x\bigr)\,dx
\;=\;
\frac{2\pi}{\sqrt{\det(-H)}},
\]
so
\begin{equation}\label{eq:local-asymp}
R_n^{\mathrm{loc}}(u)
\;=\;
e^{\Psi_u(\alpha^*,\beta^*)}\,\frac{2\pi}{\sqrt{\det(-H)}}\,\bigl(1+O(n^{-1/5})\bigr).
\end{equation}

\textit{The remaining summands.} Let $\mathcal F_n$ be the complementary lattice in $\Delta$. On the compact set $\Delta$ the continuous function $\eta$ attains the value $0$ only at $(\alpha^*,\beta^*)$, so $\eta\le-\delta$ off a fixed ball about the maximiser. In the annulus $\omega_n\le\|(\alpha,\beta)-(\alpha^*,\beta^*)\|\le\kappa$ the quadratic upper bound from $H\prec0$ gives $\eta\le -c\omega_n^2=-c\,n^{-4/5}$. Near $\partial\Delta$ the expansion \eqref{eq:rho-summand} is not uniform (the factor $-\tfrac12\log(\alpha\beta\hat\gamma)$ diverges), but the two-sided bounds $(N/e)^N\le N!\le e\sqrt{N+1}\,(N/e)^N$ and the analogous bounds for $\Gamma(N+u/2)$ with $|u-1|\le\varepsilon$ (valid for all $N\ge0$; for $N=0$, which occurs at $p=0$, $s=0$ or $\tau=0$, they read $1\le0!\le e$ and $\Gamma(u/2)$, $\Gamma(u^2/2)$ are bounded and bounded away from $0$ on the disk, while $M(0)=1$) give $|\mathrm{summand}|\le n^{C}\,e^{n\eta(\alpha,\beta)}$ on all of $\Delta$, for a constant $C=C(d,\varepsilon)$. On a fixed neighbourhood of $\partial\Delta$ one has $\eta\le-\delta/2$, so those $O(n^2)$ points contribute $O(n^{C+2}e^{-\delta n/2})$. There are $O(n^2)$ lattice points in total, and each remaining summand off that neighbourhood is $O\bigl(n^{-1}\exp(-c'n^{1/5})\bigr)$ uniformly for $|u-1|\le\varepsilon$. Hence
\begin{equation}\label{eq:far-bound}
R_n^{\mathrm{tail}}(u)
\;:=\;
\sum_{\mathcal F_n}\rho_n(\alpha,\beta;u)
\;=\;
O\bigl(n\exp(-c'n^{1/5})\bigr)
\;=\;
o(n^{-N})
\end{equation}
for every $N$. Adding the two parts,
\[
R_n(u)
\;=\;
e^{\Psi_u(\alpha^*,\beta^*)}\,\frac{2\pi}{\sqrt{\det(-H)}}\,\bigl(1+O(n^{-1/5})\bigr),
\]
uniformly on the closed disk $\overline{D}(1,\varepsilon)$.

\textit{Leading constant.} From \eqref{eq:Psi-u} at $u=1$,
\[
e^{\Psi_1(\alpha^*,\beta^*)}\,\frac{2\pi}{\sqrt{\det(-H)}}
\;=\;
\bigl(\alpha^*\beta^*\hat\gamma^*\det H\bigr)^{-1/2}.
\]
The product of the coordinates is $\alpha^*\beta^*\hat\gamma^*=4(d-2)^2(d-3)/(d^3(d-1)^3)$, and with \eqref{eq:hessian-det}
\[
\alpha^*\beta^*\hat\gamma^*\det H
\;=\;
\frac{d-2}{d-1},
\]
hence the leading constant at $u=1$ is $C(1)=\sqrt{(d-1)/(d-2)}$. The dependence on $u$ is given by the factor
\begin{align*}
e^{\Psi_u(\alpha^*,\beta^*)-\Psi_1(\alpha^*,\beta^*)}
&\;=\;
(\alpha^*+\beta^*)^{u-1}\Bigl(\frac{\sigma^*}{\beta^*}\Bigr)^{(u^2-1)/2}\\
&\;=\;
\Bigl(\frac{d-2}{d}\Bigr)^{u-1}\Bigl(\frac{d-1}{d-2}\Bigr)^{(u^2-1)/2},
\end{align*}
which equals $C(u)/C(1)$ by \eqref{eq:Cu} (an algebraic identity valid for all complex $u$). Therefore
\[
R_n(u)\;=\;C(u)\,\bigl(1+O(n^{-1/5})\bigr)
\]
uniformly on $\overline{D}(1,\varepsilon)$. Since $u\mapsto\Psi_u(\alpha^*,\beta^*)$ is holomorphic on a neighbourhood of $u=1$, no additional analytic continuation step is required.

\textit{The cubic case $d=3$.} Here $b=0$ forces $n_0=0$, and the sum is one-dimensional in $\beta\in[0,1/2]$, with $\sigma=1-\beta$, $\hat\gamma=\tau=1-2\beta$. The same Stirling computation (now $n!/(p!\,g!)$, which contributes a factor of $n^{-1/2}$) gives
\begin{equation}\label{eq:eta-d3}
\begin{aligned}
\eta(\beta)
&\;=\;
\sigma\log\sigma-\tfrac{\hat\gamma}{2}\log\hat\gamma+\beta\log2+\tfrac32\log3-\log6,\\
\Psi_u(\beta)
&\;=\;
-\tfrac12\log(2\pi)-\tfrac12\log(\beta\hat\gamma)
+(u-1)\log\beta
+\tfrac{u^2-1}{2}\log(\sigma/\beta),
\end{aligned}
\end{equation}
and $\rho_n(\beta;u)=n^{-1/2}\exp(n\eta(\beta)+\Psi_u(\beta)+o(1))$. One has $\eta'(\beta)=\log(2\hat\gamma/\sigma)$, so the unique critical point in $(0,1/2)$ is $\beta^*=1/3$, with $\eta(1/3)=0$ and $\eta''(1/3)=-9/2<0$. The one-dimensional Laplace method (neighbourhood of width $\omega_n=n^{-2/5}$, tail $\exp(-c n^{1/5})$) yields
\[
R_n(u)
\;=\;
e^{\Psi_u(1/3)}\sqrt{\frac{2\pi}{|\eta''(1/3)|}}\,\bigl(1+O(n^{-1/5})\bigr).
\]
At $u=1$,
\[
e^{\Psi_1(1/3)}\sqrt{\frac{2\pi}{9/2}}
\;=\;
\bigl(\beta^*\hat\gamma^*\,|\eta''(1/3)|\bigr)^{-1/2}
\;=\;
\bigl(\tfrac13\cdot\tfrac13\cdot\tfrac92\bigr)^{-1/2}
\;=\;
\sqrt2
\;=\;
C(1),
\]
and the $u$-dependent factor is $(\beta^*)^{u-1}(\sigma^*/\beta^*)^{(u^2-1)/2}=2^{(u-1)^2/2}(2/3)^{u-1}=C(u)/C(1)$ for $d=3$. The same formulae hold uniformly on $\overline{D}(1,\varepsilon)$.
\end{proof}

\section{Exact enumeration and simulations}\label{sec:numerics}

\paragraph{Code availability.}
All implementations and reproducibility scripts for the computations reported below are available at \url{https://github.com/ktadler/2factor-cycles}.

Small values of $n$ are treated by exhaustive enumeration: for $d=3$, every perfect-matching complement is enumerated up to $n=80$, and for $d=4$, every 2-factor is enumerated up to $n=24$. For larger $n$, each trial consists of one simple $d$-regular graph generated via \texttt{networkx.random\_regular\_graph} \cite{networkx} and one 2-factor obtained by the Karp--Sipser algorithm followed by Edmonds' blossom algorithm \cite{anastos2021,edmonds1965,karpsipser1981}, applied directly to $G$ when $d=3$ and to Tutte's gadget \cite{tutte1954} (which has $2(d-1)n$ vertices) when $d\ge 4$. While this heuristic is not guaranteed to sample a uniformly random 2-factor of $G$, its empirical mean agrees closely with the exact uniform per-graph mean throughout the range where exhaustive enumeration is feasible.

\paragraph{Exact verification.}
All identities were verified in exact arithmetic. The counts $N_c = n![x^n]L_{\ge3}(x)^c/c!$ of simple 2-regular graphs coincide with brute-force enumeration for all $3\le n\le 9$; and for $n=4$, $d=3$, the identity of Theorem~\ref{thm:moments} was verified integer-for-integer against all $M(12)=10{,}395$ pairings ($\sum_P Z_c = 11664,\,10692,\,2916,\,243$ for $c=1,\dots,4$). Monte Carlo estimates at $(n,d)=(6,3)$ and $(5,4)$ agree with the theoretical formulas within sampling error.

\paragraph{Exact enumeration on fixed graphs ($d=3$, $n\le 80$; $d=4$, $n\le 24$).}
For $d=3$, the 2-factors of a graph are the complements of its perfect matchings. For each $n\in\{8,12,16,24,32,48,64,80\}$, we sampled random simple cubic graphs, enumerated all perfect matchings by backtracking, and recorded the exact conditional mean of $X'$ on each graph; see Table~\ref{tab:typical}. For $d=4$, where no single-matching reduction is available and the expected number of 2-factors $\E Z$ is much larger (approximately $1.7\times 10^{14}$ at $n=80$, compared with $1.4\times 10^5$ for $d=3$), exhaustive backtracking was carried out for $40$ simple $4$-regular graphs at each $n\in\{8,12,16,20,24\}$, yielding per-graph sample means $1.29,1.49,1.62,1.74,1.82$ (compared with the Karp--Sipser sample means $1.33,1.39,1.64,1.83,1.80$).

\begin{table}[tbp]
\centering
\begin{tabular}{rrrrr}
\hline
$n$ & graphs & 2-factors & mean $X'$ (per-graph avg $\pm$ sd) & $\tfrac12\log n + 0.232$\\
\hline
8 & 40 & 217 & $1.33 \pm 0.16$ & 1.27\\
12 & 40 & 394 & $1.52 \pm 0.20$ & 1.47\\
16 & 40 & 729 & $1.67 \pm 0.19$ & 1.62\\
24 & 40 & 2257 & $1.88 \pm 0.25$ & 1.82\\
32 & 40 & 7161 & $2.06 \pm 0.36$ & 1.96\\
48 & 40 & 67952 & $2.26 \pm 0.26$ & 2.17\\
64 & 20 & 372512 & $2.34 \pm 0.14$ & 2.31\\
80 & 10 & 1953655 & $2.53 \pm 0.23$ & 2.42\\
\hline
\end{tabular}
\caption{Exact enumeration on fixed graphs: mean cycle count of a uniform random 2-factor per sampled cubic graph.}
\label{tab:typical}
\end{table}

\begin{figure}[tbp]
\centering
\includegraphics[width=0.72\textwidth]{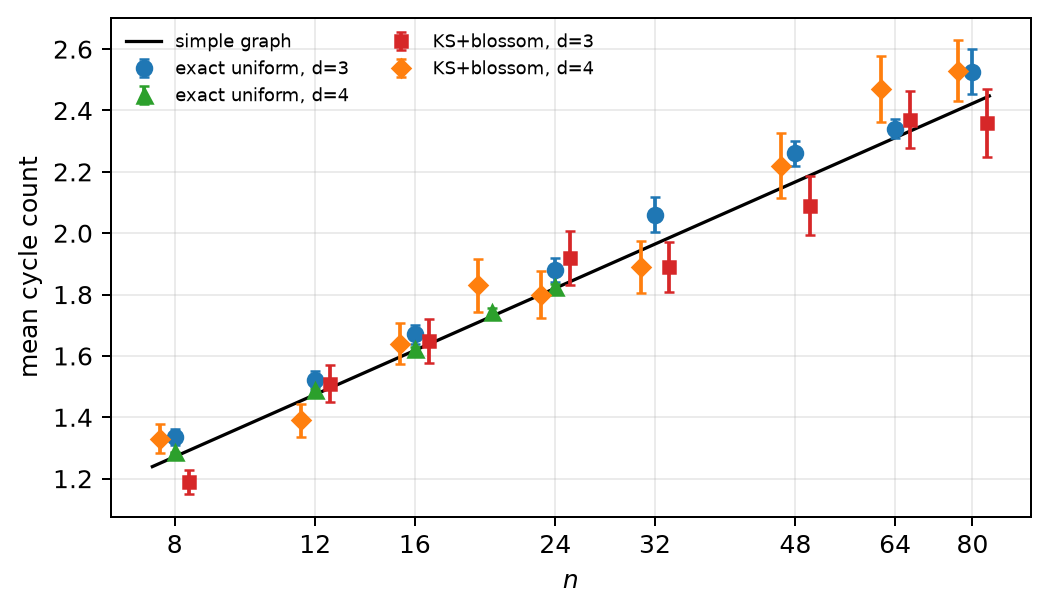}\\[0.15em]
\includegraphics[width=0.72\textwidth]{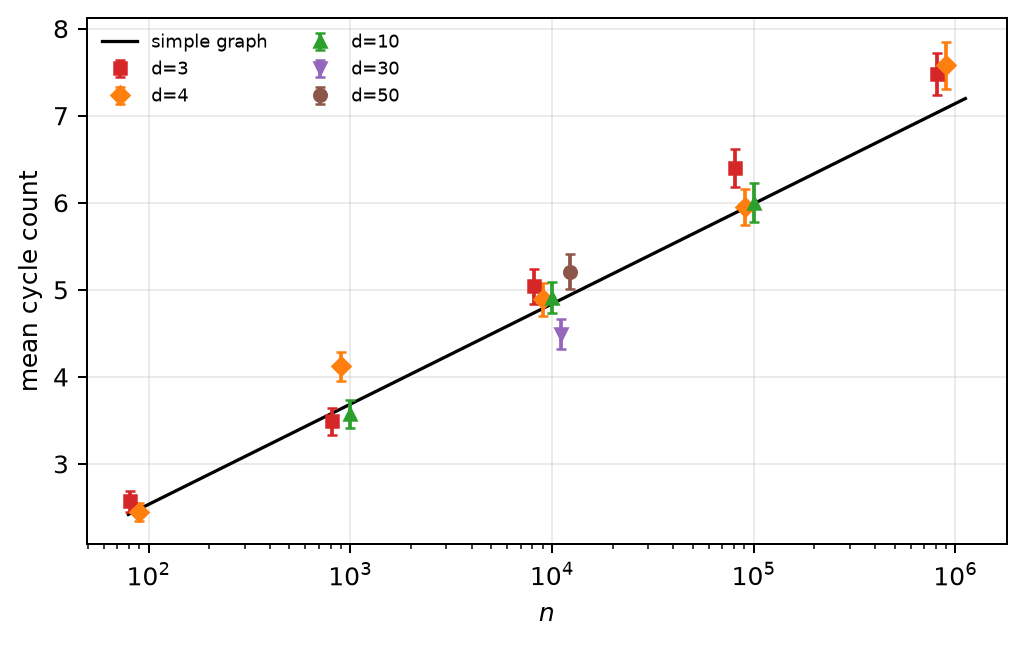}\\[0.15em]
\includegraphics[width=0.72\textwidth]{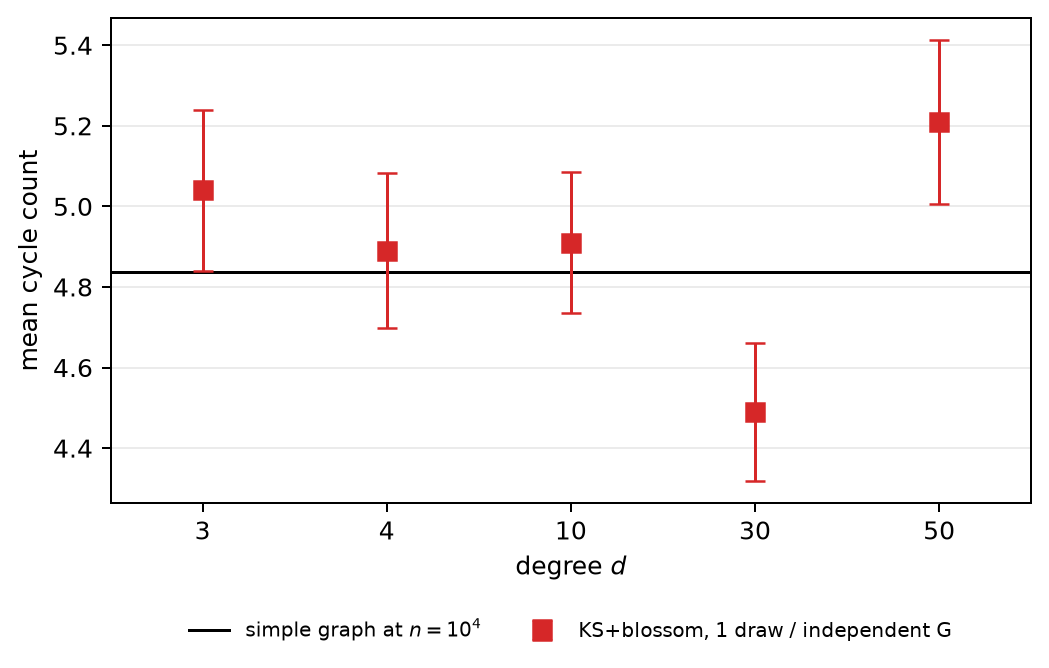}
\caption{\textbf{(a)} Comparison of exact per-graph uniform means and Karp--Sipser + blossom sample means ($N=100$) for $n\le80$ ($d=3$) and $n\le24$ ($d=4$). \textbf{(b)} Sample means ($N=100$) across powers of ten up to $n=10^6$. \textbf{(c)} Sample means at fixed $n=10^4$ for $d\in\{3,4,10,30,50\}$. Dashed line: simple-graph mean of Theorem~\ref{thm:girth}; error bars show $\pm\,\mathrm{sd}/\sqrt{N}$.}
\label{fig:sims}
\end{figure}

\begin{figure}[tbp]
\centering
\includegraphics[width=0.78\textwidth]{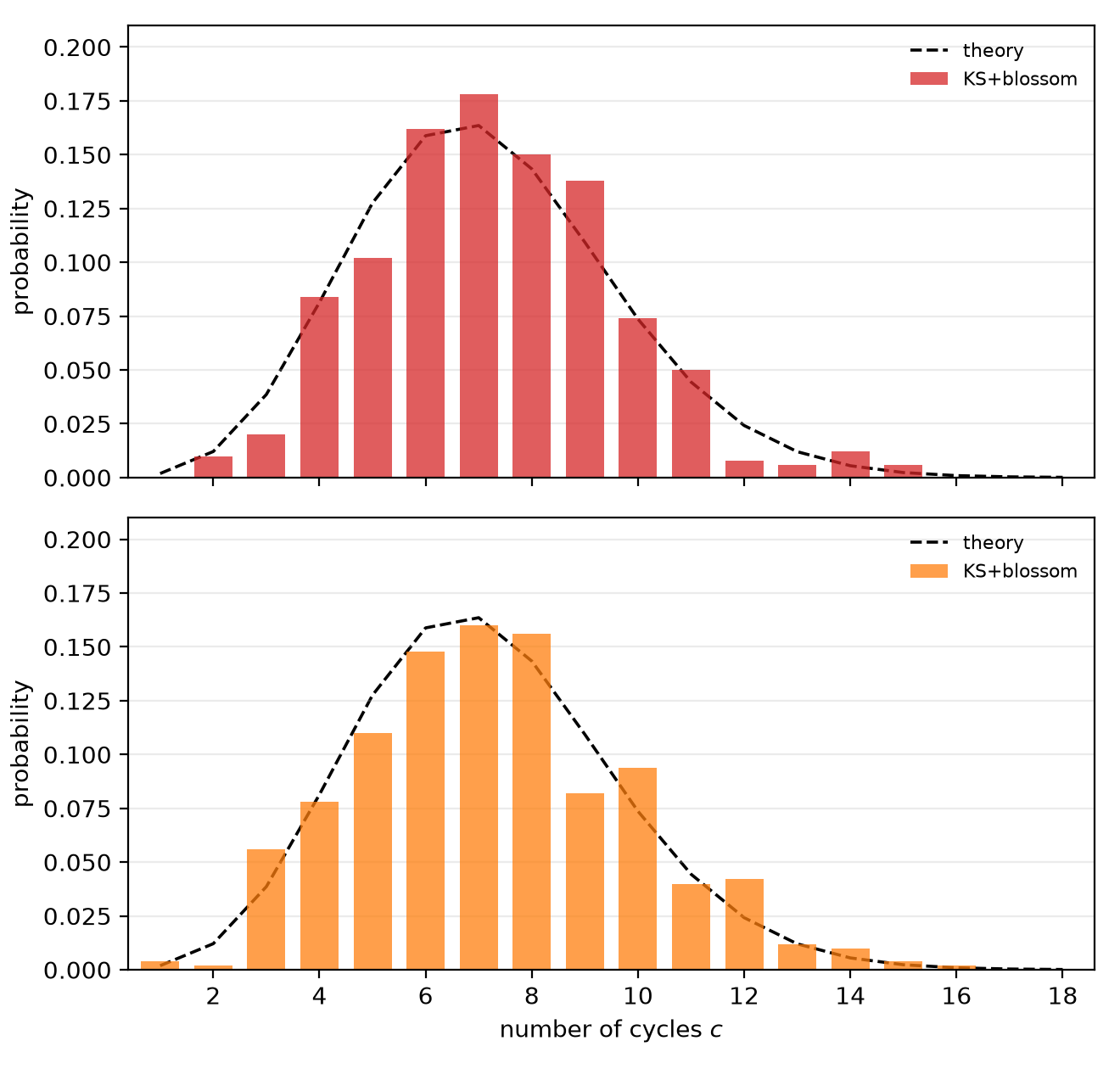}
\caption{Empirical distribution of cycle counts at $n=10^6$ ($N=500$ trials) for (a)~$d=3$ and (b)~$d=4$, compared with the theoretical distribution of Theorem~\ref{thm:girth} (dashed curve).}
\label{fig:dist}
\end{figure}

\paragraph{Simulations for larger $n$ and $d$.}
For each pair $(n,d)$ in Table~\ref{tab:simulation} and Figure~\ref{fig:sims}, we generated $N=100$ independent simple $d$-regular graphs and sampled one 2-factor per graph using the Karp--Sipser heuristic finished by Edmonds' blossom algorithm; reported uncertainties $\pm$ denote the standard error $\mathrm{sd}/\sqrt{N}$. Figure~\ref{fig:sims}(a) compares these empirical means with the exact uniform per-graph means for $n\le 80$ ($d=3$) and $n\le 24$ ($d=4$). Figure~\ref{fig:sims}(b) extends the estimates to $n\le 10^6$ for $d=3,4$ and $n\le 10^5$ for $d=10$, comparing them against the asymptotic simple-graph mean $\mu_n = \tfrac12\log n + \tfrac\gamma2 + \log 2 - \tfrac34$ of Theorem~\ref{thm:girth}. Figure~\ref{fig:sims}(c) holds $n=10^4$ fixed and varies $d\in\{3,4,10,30,50\}$, illustrating the degree-independence of the mean cycle count. Finally, Figure~\ref{fig:dist} compares the empirical distribution of the number of cycles over $N=500$ independent trials at $n=10^6$ for $d=3$ and $d=4$ with the theoretical distribution of Theorem~\ref{thm:girth}.

\begin{table}[tbp]
\centering
\small
\begin{tabular}{rrrr}
\hline
$d$ & $n$ & mean $X'$ ($\pm$ se) & simple graph\\
\hline
3 & 8 & $1.19\pm0.04$ & 1.27\\
3 & 12 & $1.51\pm0.06$ & 1.47\\
3 & 16 & $1.65\pm0.07$ & 1.62\\
3 & 24 & $1.92\pm0.09$ & 1.82\\
3 & 32 & $1.89\pm0.08$ & 1.96\\
3 & 48 & $2.09\pm0.10$ & 2.17\\
3 & 64 & $2.37\pm0.09$ & 2.31\\
3 & 80 & $2.36\pm0.11$ & 2.42\\
3 & $10^2$ & $2.57\pm0.12$ & 2.53\\
3 & $10^3$ & $3.49\pm0.16$ & 3.69\\
3 & $10^4$ & $5.04\pm0.20$ & 4.84\\
3 & $10^5$ & $6.40\pm0.22$ & 5.99\\
3 & $10^6$ & $7.48\pm0.24$ & 7.14\\
4 & 8 & $1.33\pm0.05$ & 1.27\\
4 & 12 & $1.39\pm0.05$ & 1.47\\
4 & 16 & $1.64\pm0.07$ & 1.62\\
4 & 20 & $1.83\pm0.09$ & 1.73\\
4 & 24 & $1.80\pm0.08$ & 1.82\\
4 & 32 & $1.89\pm0.09$ & 1.96\\
4 & 48 & $2.22\pm0.11$ & 2.17\\
4 & 64 & $2.47\pm0.11$ & 2.31\\
4 & 80 & $2.53\pm0.10$ & 2.42\\
4 & $10^2$ & $2.45\pm0.10$ & 2.53\\
4 & $10^3$ & $4.12\pm0.17$ & 3.69\\
4 & $10^4$ & $4.89\pm0.19$ & 4.84\\
4 & $10^5$ & $5.95\pm0.21$ & 5.99\\
4 & $10^6$ & $7.58\pm0.27$ & 7.14\\
10 & $10^3$ & $3.57\pm0.16$ & 3.69\\
10 & $10^4$ & $4.91\pm0.18$ & 4.84\\
10 & $10^5$ & $6.00\pm0.22$ & 5.99\\
30 & $10^4$ & $4.49\pm0.17$ & 4.84\\
50 & $10^4$ & $5.21\pm0.20$ & 4.84\\
\hline
\end{tabular}
\caption{Sample means from $N=100$ independent trials per $(n,d)$ pair (one 2-factor per graph), compared with the asymptotic simple-graph mean $\tfrac12\log n+\tfrac\gamma2+\log2-\tfrac34$ of Theorem~\ref{thm:girth}.}
\label{tab:simulation}
\end{table}

\FloatBarrier

\section*{Competing interests}
The authors declare none.

\section*{Data availability}
The code and the numerical comparisons reported in Appendix~\ref{sec:numerics} are openly available at \url{https://github.com/ktadler/2factor-cycles}.

\section*{Funding}
This work received no specific grant from any funding agency, commercial or not-for-profit sectors.

\end{document}